%% file: Lecture_NumLinAlg.tex
\documentclass[12pt]{article}

\input{headlines}

\title{\large\sffamily\bfseries Lecture Notes \\
\LARGE{Numerical Linear Algebra}\\ \large{Least Squares, QR and SVD}}
\author{\Large\sffamily Davoud Mirzaei\\ Uppsala University}
\date{\sffamily September 1, 2023}

\begin{document}

\setlength{\abovedisplayskip}{3pt}
\setlength{\belowdisplayskip}{3pt}
\setlength{\abovedisplayshortskip}{0pt}
\setlength{\belowdisplayshortskip}{0pt}
\maketitle

\definecolor{contcol1}{HTML}{72E094}
\definecolor{contcol2}{HTML}{24E2D6}
\definecolor{convcol1}{HTML}{C0392B}
\definecolor{convcol2}{HTML}{8E44AD}

\begin{tcolorbox}[title=Contents, fonttitle=\huge\sffamily\bfseries\selectfont,interior style={left color=contcol1!10!white,right color=contcol2!10!white},frame style={left color=contcol1!30!white,right color=contcol2!30!white},coltitle=black,top=2mm,bottom=2mm,left=2mm,right=2mm,drop fuzzy shadow,enhanced,breakable]
\makeatletter
\@starttoc{toc}
\makeatother
\end{tcolorbox}

\vspace*{10mm}

\thispagestyle{empty}
\newpage
\pagenumbering{arabic}

\input{lec2_part1}
\input{lec2_part2}

\input{lec2_part3}
\input{workouts_2}

\end{document}

%% file: headlines.tex
\usepackage[english]{babel}
\usepackage{graphicx,caption}
\usepackage{framed}
\usepackage[normalem]{ulem}
\usepackage{indentfirst}
\usepackage{amsmath,amsthm,amssymb,amsfonts}
\usepackage[italicdiff]{physics}
\usepackage[T1]{fontenc}
\usepackage{mathtools,extarrows}

\usepackage{setspace}
\usepackage{wrapfig}
\usepackage{lmodern,mathrsfs}
\usepackage[inline,shortlabels]{enumitem}
\setlist{topsep=2pt,itemsep=2pt,parsep=3pt,partopsep=2pt}

\usepackage[dvipsnames]{xcolor}
\usepackage[utf8]{inputenc}
\usepackage[a4paper, top=0.8in,bottom=0.8in, left=0.8in, right=0.8in, footskip=0.3in, includefoot]{geometry}
\usepackage[most]{tcolorbox}
\usepackage{tikz,tikz-3dplot,tikz-cd,tkz-tab,tkz-euclide,pgf,pgfplots}
\pgfplotsset{compat=newest}
\usepackage{multicol}
\usepackage[bottom,multiple]{footmisc} 

\usepackage[backend=bibtex,style=numeric]{biblatex}
\usepackage[nameinlink]{cleveref} 

\usepackage{graphicx}
 \graphicspath{{figures/}} 

\numberwithin{equation}{section}

\usepackage{framed,color}
\definecolor{shadecolor}{rgb}{0.9412,0.9725, 1}

\usepackage{caption}
\usepackage{tcolorbox}

\usepackage{multirow}
\usepackage{epstopdf}
\usepackage{listings}
\usepackage[]{mcode}
\usepackage{algorithm}
\usepackage{latexsym}
\usepackage{showidx}
\usepackage{latexsym}
\usepackage{amssymb}

\newcommand{\ignore}[1]{}
\usepackage{todonotes}

\usepackage{tikz}

\newcommand{\hide}[1]{} 

\renewcommand{\qed}{\hfill\blacksquare}

\newcommand{\C}{\mathbb{C}}
\newcommand{\R}{\mathbb{R}}

\newcommand{\fl}{\mathrm{fl}}

\newcommand{\wh}{\widehat}
\newcommand{\wt}{\widetilde}

\newcommand{\spann}{\mathrm{span}}
\newcommand{\range}{\mathrm{range}}
\newcommand{\nul}{\mathrm{null}}
\newcommand{\cond}{\mathrm{cond}}
\newcommand{\rankk}{\mathrm{rank}}
\newcommand{\diagg}{\mathrm{diag}}
\newcommand{\varr}{\mathrm{var}}
\newcommand{\covv}{\mathrm{cov}}

\newtheoremstyle{mystyle}{}{}{}{}{\sffamily\bfseries}{.}{ }{}
\newtheoremstyle{cstyle}{}{}{}{}{\sffamily\bfseries}{.}{ }{\thmnote{#3}}
\makeatletter
\renewenvironment{proof}[1][\proofname] {\par\pushQED{\qed}{\normalfont\sffamily\bfseries\topsep6\p@\@plus6\p@\relax #1\@addpunct{.} }}{\popQED\endtrivlist\@endpefalse}
\makeatother
\theoremstyle{mystyle}{\newtheorem{definition}{Definition}[section]}
\theoremstyle{mystyle}{\newtheorem{theorem}[definition]{Theorem}}
\theoremstyle{mystyle}{\newtheorem{lemma}[definition]{Lemma}}
\theoremstyle{mystyle}{}
\theoremstyle{mystyle}{}
\theoremstyle{mystyle}{}
\theoremstyle{mystyle}{}
\theoremstyle{mystyle}{\newtheorem{workout}[definition]{Workout}}
\theoremstyle{mystyle}{}
\theoremstyle{mystyle}{\newtheorem{labexercise}[definition]{Lab Exercise}}

\usepackage{amsmath}
\newtheorem{remark}{Remark}[section]
\newtheorem{example}{Example}[section]

\theoremstyle{cstyle}{}

\newtheoremstyle{warn}{}{}{}{}{\normalfont}{}{ }{}
\theoremstyle{warn}

\newcommand{\warningsign}[1]{\tikz[scale=#1,every node/.style={transform shape}]{\draw[-,line width={#1*0.8mm},red,fill=yellow,rounded corners={#1*2.5mm}] (0,0)--(1,{-sqrt(3)})--(-1,{-sqrt(3)})--cycle;
\node at (0,-1) {\fontsize{48}{60}\selectfont\bfseries!};}}

\tcolorboxenvironment{definition}{boxrule=0pt,boxsep=0pt,colback={red!10},left=8pt,right=8pt,enhanced jigsaw, borderline west={2pt}{0pt}{red},sharp corners,before skip=10pt,after skip=10pt,breakable}
\tcolorboxenvironment{proposition}{boxrule=0pt,boxsep=0pt,colback={Orange!10},left=8pt,right=8pt,enhanced jigsaw, borderline west={2pt}{0pt}{Orange},sharp corners,before skip=10pt,after skip=10pt,breakable}
\tcolorboxenvironment{theorem}{boxrule=0pt,boxsep=0pt,colback={blue!10},left=8pt,right=8pt,enhanced jigsaw, borderline west={2pt}{0pt}{blue},sharp corners,before skip=10pt,after skip=10pt,breakable}
\tcolorboxenvironment{lemma}{boxrule=0pt,boxsep=0pt,colback={blue!10},left=8pt,right=8pt,enhanced jigsaw, borderline west={2pt}{0pt}{blue},sharp corners,before skip=10pt,after skip=10pt,breakable}
\tcolorboxenvironment{corollary}{boxrule=0pt,boxsep=0pt,colback={violet!10},left=8pt,right=8pt,enhanced jigsaw, borderline west={2pt}{0pt}{violet},sharp corners,before skip=10pt,after skip=10pt,breakable}
\tcolorboxenvironment{proof}{boxrule=0pt,boxsep=0pt,blanker,borderline west={2pt}{0pt}{CadetBlue!80!white},left=8pt,right=8pt,sharp corners,before skip=10pt,after skip=10pt,breakable}
\tcolorboxenvironment{remark}{boxrule=0pt,boxsep=0pt,blanker,borderline west={2pt}{0pt}{Green},left=8pt,right=8pt,before skip=10pt,after skip=10pt,breakable}
\tcolorboxenvironment{remarks}{boxrule=0pt,boxsep=0pt,blanker,borderline west={2pt}{0pt}{Green},left=8pt,right=8pt,before skip=10pt,after skip=10pt,breakable}
\tcolorboxenvironment{example}{boxrule=0pt,boxsep=0pt,blanker,borderline west={2pt}{0pt}{Black},left=8pt,right=8pt,sharp corners,before skip=10pt,after skip=10pt,breakable}
\tcolorboxenvironment{examples}{boxrule=0pt,boxsep=0pt,blanker,borderline west={2pt}{0pt}{Black},left=8pt,right=8pt,sharp corners,before skip=10pt,after skip=10pt,breakable}
\tcolorboxenvironment{cthm}{boxrule=0pt,boxsep=0pt,colback={gray!10},left=8pt,right=8pt,enhanced jigsaw, borderline west={2pt}{0pt}{gray},sharp corners,before skip=10pt,after skip=10pt,breakable}
\tcolorboxenvironment{workout}{boxrule=0pt,boxsep=0pt,colback={Cyan!0},left=8pt,right=8pt,enhanced jigsaw, borderline west={2pt}{0pt}{Cyan},sharp corners,before skip=10pt,after skip=10pt,breakable}
\tcolorboxenvironment{miniproject}{boxrule=0pt,boxsep=0pt,colback={gray!10},left=8pt,right=8pt,enhanced jigsaw, borderline west={2pt}{0pt}{gray},sharp corners,before skip=10pt,after skip=10pt,breakable}
\tcolorboxenvironment{labexercise}{boxrule=0pt,boxsep=0pt,colback={green!10},left=8pt,right=8pt,enhanced jigsaw, borderline west={2pt}{0pt}{green},sharp corners,before skip=10pt,after skip=10pt,breakable}

\newenvironment{talign*}{\csname align*\endcsname}{\endalign}

\usepackage[explicit]{titlesec}
\titleformat{\section}{\fontsize{18}{30}\sffamily\bfseries}{\thesection}{18pt}{#1}
\titleformat{\subsection}{\fontsize{15}{18}\sffamily\bfseries}{\thesubsection}{15pt}{#1}
\titleformat{\subsubsection}{\fontsize{10}{12}\sffamily\large\bfseries}{\thesubsubsection}{8pt}{#1}

\titlespacing*{\section}{0pt}{5pt}{5pt}
\titlespacing*{\subsection}{0pt}{5pt}{5pt}
\titlespacing*{\subsubsection}{0pt}{5pt}{5pt}

\DeclareMathAlphabet\mathbfcal{OMS}{cmsy}{b}{n}
\usepackage{tkz-euclide}
\usetikzlibrary{angles,quotes}

\usetikzlibrary{shapes, calc, decorations}
\newcommand{\vsp}{\vspace{0.3cm}}

\newcommand\colortxt[2][]{\tikz[baseline=(char.base)]\node[minimum width=2em,text height=2ex, fill=yellow!100,text=black,#1](char){#2};}%

%% file: lec2_part1.tex
These lecture notes focus on some numerical linear algebra algorithms in scientific computing. We assume that students are familiar with elementary linear algebra concepts such as vector spaces, systems of equations, matrices, norms, eigenvalues, and eigenvectors. In the numerical part, we do not pursue Gaussian elimination and other LU factorization algorithms for square systems. Instead, we mainly focus on overdetermined systems, least squares solutions, orthogonal factorizations, and some applications to data analysis and other areas. The reference textbooks \cite{Datta:2010}, \cite{Heath:2018}, and \cite{Trefethen-Bau:1997} are our main sources in this lecture.

\section{Least squares problem}

Let \( A \in \mathbb{R}^{m \times n} \), where \( m > n \). The system
\[
Ax  = b
\]
for a given vector \( b \in \mathbb{R}^m \) and solution \( x  \in \mathbb{R}^n \) is termed {\bf overdetermined} because it contains more equations than unknowns. This system essentially asks whether \( b \) can be expressed as a linear combination of the columns of \( A \):
$$
b = x_1a_{\cdot 1} + x_2a_{\cdot 2} + \cdots + x_na_{\cdot n}.
$$
For the square system (the case \( m = n \)), the answer is ``yes'' provided that the column vectors \( \{a_{\cdot  1}, a_{\cdot  2}, \ldots, a_{\cdot  n}\} \) are linearly independent (or equivalently for a nonsingular matrix \( A \)). However, for overdetermined systems (\( m > n \)), the answer is usually ``no'' unless \( b \) happens to lie in the span of \( \{a_{\cdot  1}, a_{\cdot  2}, \ldots, a_{\cdot  n}\} \) (often denoted \( \text{span}(A) \) or \( \text{range}(A) \)), which is highly unlikely in most applications. Therefore, in general, such a system has no solution.

To illustrate this, consider the case where \( m = 3 \) and \( n = 2 \). In this scenario, \( a_{\cdot  1} \) and \( a_{\cdot  2} \) represent two vectors in \( \mathbb{R}^3 \). If \( a_{\cdot  1} \) and \( a_{\cdot  2} \) are linearly independent, then their span forms a plane (a \( 2 \)-dimensional subspace) in \( \mathbb{R}^3 \). The system \( Ax  = b \) has a solution if \( b \) lies in that plane; otherwise, the system has no solution. The probability of a vector \( b \in \mathbb{R}^3 \) lying in a plane is zero.

In such situations, one obvious alternative to ``solving the linear system exactly'' is to minimize the residual vector
\[
r = b - (x_1a_{\cdot  1} + x_2a_{\cdot  2} + \cdots + x_na_{\cdot  n}) = b - Ax .
\]
The solution to the problem depends on how we measure the length of the residual vector. It is preferred to use the \( 2 \)-norm, although any norm could be used. The \( 2 \)-norm is induced by the inner product, thus is related to the notion of orthogonality, and it is smooth and strictly convex. These properties make the theory and computation with this norm much easier than with other norms. With the use of the \( 2 \)-norm, the solution is the vector \( x  \) that minimizes the sum of squares of differences between the components of \( b \) and \( Ax  \). This method is known as {\bf least squares}. In the least squares method, we seek to find an optimal vector that solves the minimization problem:
\begin{equation}\label{leastsquaredef}
\min_{x  \in \mathbb{R}^n} \|Ax  - b\|_2. \vsp
\end{equation}
As we pointed out, the solution \( x  \) of this problem (which always exists) may not exactly satisfy \( Ax  = b \). To reflect the lack of exact equality, we may write the linear least squares problem as
\[
Ax  \cong b,
\]
and approximation is understood in the least square sense.
\vsp

\begin{example}[\cite{Heath:2018}] \label{ex:height_hills}
A surveyor tries to measure the heights of three hills. Sighting first, his/her initial measurements are
$x_1 = 1237$ ft, $x_2=1941$ ft, and $x_3=2417$ ft. To confirm these measurements, the surveyor climbs to the top of the first hill and measures the heights of the second and third hills above the first and obtain $x_2-x_1=711$ ft and $x_3-x_1=1177$ ft. Then he/she climbs to the top of the second hill and measures $x_3-x_2 = 475$ ft. It is obvious that there exists an inconsistency with measurements. These
can be written in an overdetermined linear system of equations:
\begin{equation*}
  Ax  = \begin{bmatrix}1&0&0\\0&1&0\\0&0&1\\-1&1&0\\-1&0&1\\0&-1&1  \end{bmatrix}
  \begin{bmatrix}
    x_1 \\
    x_2\\
    x_3
  \end{bmatrix}
  \cong\begin{bmatrix}
     1237 \\
     1941 \\
     2417\\
     711\\
     1177\\
     475
   \end{bmatrix} = b.
\end{equation*}
The least squares solution of this system (as we will see soon) is $x = [1236, 1943, 2416]^T$ which differs slightly from the initial measurements. The last three observations helped to obtain a better measurement.
\end{example}
\vsp
\begin{example}
Data fitting (or curve fitting) is a procedure for finding the curve of best fit to a given set of data points.
Trying to find the best curve by minimizing the sum of the squares of the residuals of the points from the curve
leads to a linear least squares problem of the form \eqref{leastsquaredef}.
Given data $(t_k,y_k)$, $k=1,2,\ldots,m$, we wish to find a function 
$p\in \spann\{\phi_1(t),\phi_2(t),\ldots,\phi_n(t)\}$ such that $p$ is the best fit to data values $y_k$ in the sense that
$$
\sum_{k=1}^{n}(y_k-p(t_k))^2 \rightarrow \min.
\vsp
$$
By expanding $p$ in terms of basis functions $\phi_j$ with coefficients $c_j$, i.e.,
$$
p(t) = c_1\phi_1(t) + c_2\phi_2(t) + \cdots + c_n\phi_n(t),
$$
the problem is equivalent with finding a vector $c=(c_1,\ldots,c_n)^T\in \R^n$ that solves the minimization problem
\begin{equation*}
  \min_{c\in\R^n} \sum_{k=1}^{m} \big( y_k - (c_1\phi_1(t_k) + \cdots + c_n\phi_n(t_k)) \big)^2.
\end{equation*}
If we define the $m\times n$ matrix $A$ with entries $a_{kj} = \phi_j(t_k)$ and $m$-vector $b=(y_1,\ldots,y_m)^T$, then the above data-fitting problem takes the form
$$
\min_{c\in\R^n}\|Ac-b\|_2.
$$
In the case where the approximation space is the space of polynomials (with basis $\{1,t,\ldots,t^{n-1}\}$, or any other basis), the problem is known as {\em polynomial curve fitting}.
A schematic of a cubic curve fitting ($n=4$) is illustrated in Figure \ref{fig:curve_fit}.

\begin{center}
  \includegraphics[scale=1]{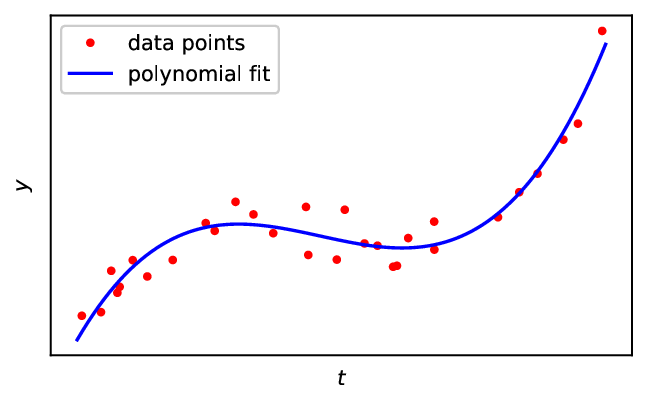}
  \captionof{figure}{A least squares polynomial fit to a given data set}\label{fig:curve_fit}
\end{center}

The special case with basis $\{1,t\}$ is refereed to as {\em linear curve fitting} or {\em linear regression}.
As an example, to find a cubic fit $p_3(t)=c_1+c_2t+c_3t^2+c_4t^3$ to six points $(t_1,y_1),\ldots,(t_6,y_6)$ the matrix $A$ is $6\times 4$ and the problem has the form
\begin{equation*}
  Ac = \begin{bmatrix}1&t_1&t_1^2 & t_1^3\\1&t_2&t_2^2& t_2^3\\1&t_3&t_3^2& t_3^3\\1&t_4&t_4^2& t_4^3\\1&t_5&t_5^2 & t_5^3\\
  1&t_6&t_6^2 & t_6^3
   \end{bmatrix}
  \begin{bmatrix}
    c_1 \\
    c_2\\
    c_3 \\
    c_4
  \end{bmatrix}
  \cong\begin{bmatrix}
     y_1 \\
     y_2 \\
     y_3\\
     y_4\\
     y_5 \\
     y_6
   \end{bmatrix} = b.\vsp 
\end{equation*}
The matrix $A$ here is a {\em Vandemonde} matrix which will be ineffectively ill-conditioned for higher order polynomial curve fittings.

\end{example}
\vsp

\subsection{Existence and uniqueness}
Assume that $y=Ax $ so $y\in \range(A)$. The function $f(y)=\|b-y\|_2$ is continuous and coercive on $\R^m$, so it has at least a
minimum on the closed and unbounded set $\mathrm{range}(A)$.  Moreover, $f$ is strictly
convex on the convex set $\mathrm{range}(A)$, so the minimum vector $y$ is unique. It does not mean that
the solution $x $ of the least square problem \eqref{leastsquaredef} is unique in general.
Suppose $x $ and $\widetilde{x }$ are two solutions for the least square problem and $z=x -\widetilde x $. Then
$Az = 0$ because $Ax =A\widetilde x $. If the columns of $A$ are linearly independent
then $z = 0$ and $x =\widetilde{x }$.  We conclude that if
$A$ has full column rank, i.e., $\mathrm{rank}(A) = n$  then the solution of least squares problem is unique
(the inverse is also true).
If $\mathrm{rank}(A) < n$, then A is said to be {\em rank-deficient}.
\vsp

\subsection{Normal equations}

 To find the solution of the least squares problem \eqref{leastsquaredef}, we define the \( n \)-variate function \( f:\mathbb{R}^n \to \mathbb{R} \) by
\[
f(x ) := \|Ax -b\|_2^2 = (Ax -b)^T(Ax -b) = x ^TA^TAx -2x ^TA^Tb + b^Tb
\]
and aim to minimize it on \( \mathbb{R}^n \). The function \( f \) is quadratic, and a necessary condition for a minimizer \( x  \) is that \( \nabla f(x ) = 0 \), i.e., \( 2A^TAx  -2A^Tb = 0 \). Therefore, any minimizer of \( f \) should satisfy
\begin{equation}\label{LS-normaleq}
   A^TAx  = A^Tb.
 \end{equation}
The {\em Hessian} of \( f \) is \( 2A^TA \), which is semi-positive definite in general. However, if the columns of \( A \) are linearly independent (meaning \( A \) has full rank), then \( A^TA \) is positive definite (you should prove this!). In this case, \eqref{LS-normaleq} becomes a sufficient condition as well. The linear system \eqref{LS-normaleq} is known as the system of {\bf normal equations} and  suggests a method for solving the least squares problem. However, it is advisable to avoid it in favor of other computationally more stable algorithms, as we will describe later.

\begin{figure}[!th]
\begin{center}
\begin{tikzpicture}
    \tkzDefPoint(0,0){A}
    \tkzDefPoint(6,0){B}
    \tkzDefPoint(10,2.5){C}
    \tkzDefPoint(4,2.5){D}
    \tkzDefPoint(5,5){E}
    \tkzDefPoint(5,1){F}
    \tkzDrawPolygon[thick,fill=gray!10](A,B,C,D)
    \draw[->,color=blue,thick] (0,0)--(5,5)node[above] {$b$};
    \draw[->,color=red,thick] (0,0)--(5,1)node[black, right] {$Ax$};
    \draw[dashed,color=red,thick] (5,5)--(5,1)node[right] {};
    \tkzLabelSegment[above,pos=.4,sloped](E,F){$b-Ax $}
    \draw[->,color=black,thick] (0,0)--(5,0)node[black, below] {$a_{\cdot  1}$};
    \draw[->,color=black,thick] (0,0)--(1.6,1)node[black, right] {$~a_{\cdot  2}$};

    \draw[color=red,thick] (4.75,.98)--(4.75,1.25) node[black, right] {};
    \draw[color=red,thick] (4.75,1.25)--(5,1.30) node[black, right] {};

\end{tikzpicture}
\end{center}
\caption{Geometric interpretation of the least squares problem.}\label{fig:bestapp}
\end{figure}
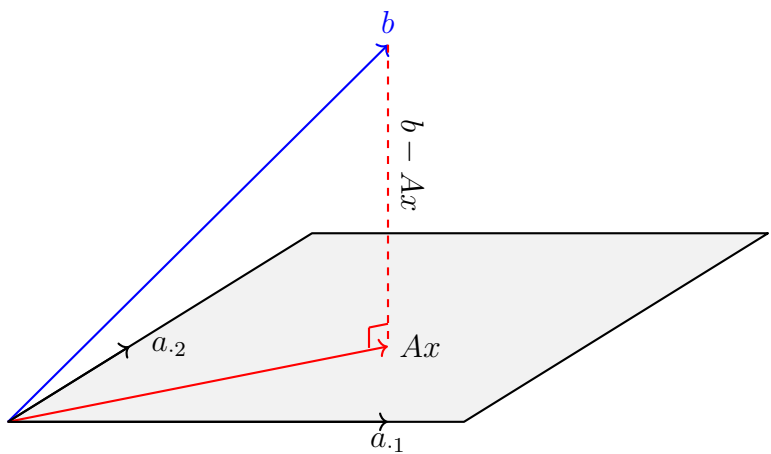

 Another approach, equivalent with one we derived above, is as below.
 The vector $Ax =x_1a_{\cdot 1}+\cdots+x_na_{\cdot  n}$ out of subspace $\range(A)$ closest to $b$ in the Euclidean norm occurs when the residual vector $b-Ax $ is perpendicular to $\range(A)$. See Figure \ref{fig:bestapp}. Thus, the inner product of $b-Ax $ and any column of $A$ should be zero, equivalently
 $$
 A^T(b-Ax )=0
 $$
which is the same system of normal equations \eqref{LS-normaleq}.

\begin{example}\label{ex:height_hills2}
Returning to Example \ref{ex:height_hills}, since
$A$ has full rank, we can obtain the least squares solution by solving the normal system $A^TA=A^Tb$. We have
$$
A^TA = \begin{bmatrix}  3&-1&-1\\ -1&3&-1\\-1&-1&3\end{bmatrix}, \quad A^Tb =
\begin{bmatrix}
-651\\
2177\\
4069
 \end{bmatrix}
$$
Since matrix $A^TA$ is positive definite, the solution of the normal system can be obtained by the Cholesky factorization $A^TA = LL^T$ where $L$ is a lower triangular matrix. We therefore need to solve two triangular systems $Lz = A^Tb$ (lower triangular) and $L^Tx  = z$ (upper triangular) using forward and backward substitutions for final solution $x $. 
In Python write
\begin{shaded}
\vspace*{-3mm}
\begin{verbatim}
import numpy as np
import scipy as sp
A = np.array([[1,0,0],[0,1,0],[0,0,1],[-1,1,0],[-1,0,1],[0,-1,1]])
b = np.array([1237,1941,2417,711,1177,475])
L = np.linalg.cholesky(A.T@A)   # Cholesky factorization
z = sp.linalg.solve_triangular(L, A.T@b, lower = True) # forward
x = sp.linalg.solve_triangular(L.T, z, lower = False)  # backward
print('Hill heights =', x)
\end{verbatim}
\vspace*{-3mm}
\end{shaded}
The final solution will be
$x = [1236, 1943, 2416]^T$.
From a computational point of view, the normal equation is effective for solving this small-sized least squares system. However, in practical scenarios and for larger matrix sizes, solving through the normal equation is strongly discouraged, as we will discuss it later.
\end{example}
\vsp

\subsection{Conditioning of the least squares problem}
The conditioning of linear system $Ax =b$ for a square and nonsingular matrix $A$ is measured by the condition number
$$
\cond(A) = \|A\|\|A^{-1}\|.
$$
For the overdetermined system $Ax \cong b$, however, the inverse of $A$ can not be defined in the conventional sense, but it is possible to define a {\bf pseudoinverse} or {\bf Moore–Penrose inverse}, denoted
by $A^{+}$, that behaves like an inverse in many respects\footnote{Such definition of inverse matrix was independently described by E. H. Moore in 1920, Arne Bjerhammar in 1951, and Roger Penrose in 1955. Earlier, Erik Ivar Fredholm had introduced the concept of a pseudoinverse of integral operators in 1903.}.
The pseudoinverse $A^{+}$ exists for any matrix $A$, in particular, when $A$ has linearly independent columns then
$$
A^{+} = (A^TA)^{-1}A^T
$$
which is indeed a {\em left inverse} because $A^{+}A=I$. On the other hand, $P:=AA^{+}$ is an orthogonal projector onto $\range(A)$, so that the solution of the least squares problem can be written as

$$
x  = A^{+}b.
$$
In Section \ref{sect:svd}, we explore how the pseudoinverse can be computed for any matrix $A$ using the SVD.
Now, for a matrix $A\in\R^{m\times n}$ with $m\geqslant n$ and $\rank(A)=n$, the condition number is defined as
\begin{equation*}
  \cond_2(A):= \|A\|_2\|A^{+}\|_2.
\end{equation*}
This definition remains valid even if $\text{rank}(A) < n$, provided that we can compute $A^{+}$. 


While the conditioning of a square system depends solely on $A$, the conditioning of a least squares problem $Ax  \cong b$ relies on both the coefficient matrix $A$ and the right-hand side vector $b$. In fact, the closeness of $b$ to $\text{range}(A)$ will affect the conditioning. From 
Figure \ref{fig:bestapp} we observe 
$$
\frac{\|Ax \|_2}{\|b\|_2} = \cos\theta, 
$$
where $\theta$ is the angle between $Ax $ and $b$.
To measure the sensitivity of the least squares solution to input perturbations, we consider perturbations in $b$ and $A$ separately.
For the perturbed vector $b+\delta$ the least squares solution $x +\varepsilon$ is given by the normal equation as
$A^TA(x +\varepsilon) = A^T(b+\delta)$.
Combining with $A^TAx  = A^Tb$ yields $A^TA\varepsilon = A^T\delta$ or $\varepsilon = A^{+}\delta$. This gives
$$
\|\varepsilon\|_2\leqslant \|A^{+}\|_2 \|\delta\|_2.
$$
Dividing both sides by $\|x \|$ then yields
\begin{align*}
  \frac{\|\varepsilon \|_2}{\|x \|_2} \leqslant & \;\|A^{+}\|_2 \frac{\|\delta \|_2}{\|x \|_2} \\
  =& \;\cond(A) \frac{\|b\|_2}{\|A\|_2\|x \|_2}\frac{\|\delta \|_2}{\|b\|_2}  \\
  \leqslant &\; \cond(A) \frac{\|b\|_2}{\|Ax \|_2}\frac{\|\delta \|_2}{\|b\|_2}  \\
  =&\; \cond(A) \frac{1}{\cos\theta}\frac{\|\delta \|_2}{\|b\|_2}.
\end{align*}
We observe that the condition number for the least squares solution $x $ with respect to perturbations
in $b$ depends on both $\cond(A)$ and the angle $\theta$ between $b$ and $Ax $. In particular, the condition number is approximately $\cond(A)$ when the
residual is small ($\cos\theta \approx 1$), however, it can be arbitrarily
worse than $\cond(A)$ when the residual is large ($\cos\theta \approx 0$).
\vsp

\begin{workout}
Assume that the perturbed least squares solution $x  + \varepsilon $ is obtained by perturbing the input matrix $A + E$, while the right-hand side vector $b$ remains unchanged. Show that
\begin{equation}\label{cond_ls_AE}
\frac{\|\varepsilon \|_2}{\|x \|_2} \leqslant \left([\cond(A)]^2\tan\theta + \cond(A) \right)\epsilon_A + \mathcal O(\epsilon_A^2)
\vsp
\end{equation}
where $\epsilon_A = \frac{\|E\|_2}{\|A\|_2}$.
Conclude that the condition number is approximately $\text{cond}(A)$ when the residual is small. However, the condition number is squared for a moderate residual, and it becomes arbitrarily large when the residual is large.
\end{workout}
\vsp

The most striking feature of \eqref{cond_ls_AE} is that it depends on the square of $\text{cond}(A)$. This implies that even if $A$ is only mildly ill-conditioned, a small perturbation in $A$ can cause a large change in $x $. An exception occurs in problems where the least squares solution fits the data very well, i.e., $\tan\theta\approx 0$, causing the factor $[\text{cond}(A)]^2$ to cancel out.
\vsp

\begin{example}
Let us again consider the height measurements of hills in Example \ref{ex:height_hills} with the least squares solution $x  = [1236, 1943, 2416]^T$. The pseudoinverse of $A$ is given by
$$
A^+ = (A^TA)^{-1}A^T=\frac{1}{4}
\begin{bmatrix}
  2 & 1 & 1 & -1 & -1 & 0 \\
 1 & 2 & 2 & 1 & 0 & -1 \\
 1 & 1 & 2 & 0 & 1 & 1
 \end{bmatrix}.
$$
We have $\|A\|_2= 2$ and $\|A^+\|_2=1$, so that
$$
\cond(A) = \|A\|_2\|A^+\|_2 = 2.
$$
On the other side, the ratio is computed as
$$
\cos\theta = \frac{\|Ax \|_2}{\|b\|_2} \doteq \frac{3640.8761}{3640.8809} \doteq 0.99999868,
$$
Hence, the angle between $Ax $ and $b$ is the small value $\theta = 0.001625$, indicating that the norm of the residual $r=b-Ax $ is very small. Given the small condition number and the angle $\theta$, we can conclude that this particular least squares problem is well-conditioned.
\end{example}
\vsp 
\begin{example}
Consider a least squares problem with coefficient matrix
$$
A = \begin{bmatrix}
      1 & 1 \\
      \epsilon & -\epsilon \\
      0 & 0
    \end{bmatrix}
$$
and right-hand side vector $b = [1,0,\epsilon]^T$, where $\epsilon>0$ is a small parameter. The pseudoinverse of $A$ is computed as
$$
A^+ = \frac{1}{2}\begin{bmatrix}
      1 & \frac{1}{\epsilon} & 0 \\
      1 & -\frac{1}{\epsilon} & 0
    \end{bmatrix}.
    \vsp
$$
The condition number of $A$ is
$$
\cond(A) = \|A\|_2\|A^+\|_2 = \sqrt{2} \frac{1}{\epsilon \sqrt 2} = \frac{1}{\epsilon},
$$
and the least squares solution is given by $x =A^+b = [1/2,1/2]^T$. Assume a tiny perturbation 
$$
E = \begin{bmatrix}
      0 & 0 \\
      0 & 0 \\
      -\epsilon & \epsilon
    \end{bmatrix}
$$
for matrix $A$. The perturbed solution then is obtained as
$$
x +\varepsilon  = (A+E)^+ b = \frac{1}{2}\begin{bmatrix}
      1 & \frac{1}{2\epsilon} & -\frac{1}{2\epsilon} \\
      1 & -\frac{1}{2\epsilon} & \frac{1}{2\epsilon}
    \end{bmatrix}
    \begin{bmatrix}
      1 \\
      0\\
      \epsilon
    \end{bmatrix}
=\begin{bmatrix}
   \frac{1}{4} \\
   \frac{3}{4}
 \end{bmatrix},
$$
which shows that $\varepsilon  = [-1/4,1/4]$ and $\frac{\|\varepsilon\|_2}{\|x \|_2}= 1/2$. We can observe this large error from stability 
 bound \eqref{cond_ls_AE}.
We have $\frac{\|E\|_2}{\|A\|_2}=\frac{\epsilon\sqrt{2}}{\sqrt{2}}=\epsilon$, and $\theta = \cos^{-1}\frac{\|Ax \|_2}{\|b\|_2}= \cos^{-1}\frac{1}{\sqrt{1+\epsilon^2}}\approx 0$ for small values of $\epsilon$. Therefore, the term with the squared condition number in the right-hand side of \eqref{cond_ls_AE} becomes negligible. Consequently, the bound on the output perturbation $\frac{\|{\varepsilon} \|_2}{\|x \|_2}$ is solely determined by $\cond(A)\frac{\|E\|_2}{\|A\|_2} = 1$, which aligns with the exact value of $\frac{\|{\varepsilon} \|_2}{\|x \|_2}$.

Now, let us change the right-hand side to $b = [1,0,1]^T$. For this case, we have
$$
x  = A^{+}b =
\begin{bmatrix}
\frac{1}{2} \\
\frac{1}{2}
\end{bmatrix},
\quad
x +\varepsilon  = (A+E)^+ b = \begin{bmatrix}
                           \frac{1}{2}-\frac{1}{4\epsilon} \\
                           \frac{1}{2}+\frac{1}{4\epsilon}
                         \end{bmatrix},
$$
and $\frac{\|{\varepsilon} \|_2}{\|x \|_2}=\frac{1}{2\epsilon}$. The relative perturbation in the solution is approximately $[\cond(A)]^2\frac{\|E\|_2}{\|A\|_2}$. In this case, $\theta = \cos^{-1}\frac{1}{\sqrt 2}=\frac{\pi}{4}$, and $\tan\theta = 1$, indicating that the condition-squared term in the perturbation bound is not suppressed. As a result, the solution is highly sensitive to perturbations.
\end{example}
\vsp

\input{projectors}

\section{Orthogonality and QR factorization}
From several points of view, it is advantageous to use orthogonal vectors as
basis vectors in a vector space. As an application, we will obtain an stable algorithm for solving
the least squares problem when a specific orthogonal basis is obtained for the subspace $\range(A)$.

We remind that two nonzero vectors $x$ and $y$ are called orthogonal if $x^T y = 0$.

\begin{workout}
If vectors $q_j$, $j=1,2,\ldots,m$ are mutually orthogonal, i.e. $q_j^Tq_k=0$ for $j\neq k$, then they
 are linearly independent.
\end{workout}
If the set of orthogonal vectors $q_j\in\R^m,\, j = 1, 2, \ldots , m$, be normalized by $\|q_j\|_2=1$
then they are called {\bf orthonormal}, and form an orthonormal basis for $\R^m$.
\begin{definition}
A square matrix whose columns are orthonormal is called an orthogonal matrix.
\end{definition}
It is clear that an orthogonal matrix $Q$ satisfies $Q^TQ = I$, it is full rank, and its inverse is equal to $Q^{-1} = Q^T$.
The rows of an orthogonal matrix are orthogonal, i.e., $QQ^T = I$.
It is not difficult to show that the product of two orthogonal matrices is orthogonal.

One of the most important properties of orthogonal matrices is that they
preserve the length (Euclidian norm) of a vector:
$$
\|Qx\|_2^2 = (Qx)^TQx=x^TQ^TQx = x^Tx = \|x\|_2^2.
$$
This means that $Q$ rotates $x$ but does not change its length.
In numerical point of view, this norm-preserving property means that orthogonal matrices do not amplify
errors.

\begin{workout}
Show that orthogonal matrices preserve the $2$-norm and the Frobenius norm of matrices.
\end{workout}

Here we introduce two classes of elementary orthogonal matrices that will be used in the sequel to transfer the columns of an arbitrary matrix
$A$ into an set of orthonormal bases.
\vsp

\subsection{Householder transformations}
In this section, we consider Householder transformations as a class of orthogonal matrices, which are particularly useful for performing reflections and orthogonal projections in matrix computations. The we will use them to compute the QR factorization of a matrix. 
\begin{definition}
A matrix of the form
\begin{equation}\label{def:housh}
  H = I -\frac{2}{u^Tu}uu^T
\end{equation}
where $u$ is a non-zero vector in $\R^n$ is called a {Householder matrix} or a {\bf Householder transformation}\footnote{After the celebrated numerical analyst Alston Scott Householder (5 May 1904 – 4 July 1993)}. The vector $u$ determining the Householder matrix
$H$ is called the Householder vector.
\end{definition}

\begin{center}
\begin{center}
\begin{tikzpicture}[scale = 0.8]
    \tkzDefPoint(0,0){A}
    \tkzDefPoint(7,0){B}
    \tkzDefPoint(10,2){C}
    \tkzDefPoint(3,2){D}
    \tkzDefPoint(4,5){E}
    \tkzDefPoint(4,3){U}
    \tkzDefPoint(4,1){O}
    \tkzDefPoint(6.5,5){X}
    \tkzDefPoint(6.5,-3){W}
    \tkzDefPoint(6.5,1){P}
    \tkzDefPoint(6.5,0){Q}
    \tkzDefPoint(4.6,0){Z}
    \tkzDefPoint(4,-3){Y}
    \tkzDrawPolygon[thick,fill=gray!10](A,B,C,D)
    \node[color=black] at (1.5,.5) {$U$};
    \tkzDrawSegments[red,dashed,->](O,E)
    \tkzDrawSegments[red,->,thick](O,U)
    \tkzDrawSegments[blue,->,thick](O,X)
    \tkzLabelSegment[right,pos=.9](O,U){$u$}
    \tkzLabelSegment[above,pos=1](O,X){$x$}
    \tkzDrawSegments[red](X,P)
    \tkzDrawSegments[red,dashed](X,E)
    \tkzDrawSegments[red,dashed](P,Q)
    \tkzDrawSegments[red,->](Q,W)
    \tkzDrawSegments[blue,dashed,thick](O,Z)
    \tkzDrawSegments[blue,->,thick](Z,W)
    \tkzDrawSegments[red,dashed](O,Y)
    \tkzLabelSegment[below,pos=.2,sloped](E,O){$u(u^Tx)$}
    \tkzLabelSegment[above,pos=.5,sloped](Q,W){$-2u(u^Tx)$}
    \tkzLabelSegment[right,pos=1](Z,W){$w=x-2uu^Tx=Hx$}
\end{tikzpicture}
\end{center}
\captionof{figure}{Geometric interpretation of Householder transformation}\label{fig:househ}
\end{center}

\vsp
See Figure \ref{fig:househ} for a geometric interpretation of a Householder transformation. In the figure the vector $u$ is assumed to be a normal vector ($u^Tu=1$) of plane $U$. The vector $w=Hx$ is
the reflection of $x$ with respect to plane $U$. The plane acts as a {\bf mirror}; the reason why the Householder transformation is also known as {\em elementary reflector}.

The following properties can be simply proved for a Householder transformation.

\begin{theorem}\label{thm:househ_peroperties}
Let $H=I-2uu^T/u^Tu$ be a Householder matrix with vector $u\in\R^n$. Then
\begin{enumerate}
  \item $H$ is symmetric and orthogonal.
  \item $H^2 = I$
  \item $Hu=-u$
  \item $Hv=v$ if $u^Tv=0$
  \item If $x,y\in\R^n$ are such that $x\neq y$ and $\|x\|_2=\|y\|_2$, and $u$ is chosen parallel to $x-y$ then $Hx=y$.
\end{enumerate}
\end{theorem}
\proof
We only prove item (5) because other items are easy to prove. Let $u = c(x-y)$ for a constant $c\neq 0$, and write $x = \frac{1}{2}(x+y)+\frac{1}{2}(x-y)$. Then
$$
Hx = \frac{1}{2}H(x+y) + \frac{1}{2}H(x-y).
$$
By using property (3) we have $H(x-y)=-(x-y)$. On the other hand $(x+y)$ is orthogonal to $x-y$ because
$(x+y)^T(x-y)=x^Tx-x^Ty + y^Tx - y^Ty = \|x\|_2-\|y\|_2=0$. Thus by property (4) we have $H(x+y)=x+y$.
All these together give $Hx=y$.
$\qed$
\vsp

\begin{remark}
Properties (3) and (4) show that $H$ has $n-1$ eigenvalues $1$ and a simple eigenvalue $-1$ (corresponding to eigenvector $u$).
\end{remark}
\vsp 

\begin{theorem}\label{thm:househ_zero}
Given a nonzero vector $x\neq e_{\cdot 1}:=[1,0,\ldots,0]^T$, the Householder matrix $H$ define by $u = x \pm \|x\|_2e_{\cdot 1}$ is such that $Hx = \mp \|x\|_2e_{\cdot 1}$.
(Take care of signs $\pm$ and $\mp$).
\end{theorem}
\proof
This is a simple consequence of item (5) of Theorem \ref{thm:househ_peroperties} by letting $y = \mp \|x\|_2e_{\cdot 1}$.
$\qed$
Here is an illustration ($\alpha=\mp \|x\|_2$):

$$
x =
\begin{pmatrix}
  x_1 \\
  x_2 \\
  x_3 \\
  \vdots \\
  x_n
\end{pmatrix}
\quad \Longrightarrow  \quad
Hx =
\begin{pmatrix}
  \alpha \\
  0 \\
  0 \\
  \vdots \\
  0
\end{pmatrix}
$$
\vsp 

Given the vector $x$, we find a reflection matrix $H$ (the {mirror}) such that $Hx = \alpha e_{\cdot 1}$ (reflect $x$ on $x_1$-axis). See Figure \ref{fig:mirror} below.

\begin{center}
\begin{center}
\begin{tikzpicture}[scale = 0.7]
    \tkzDefPoint(0,0){O}
    \tkzDefPoint(0,5){Z}
    \tkzDefPoint(4,-1){Y}
    \tkzDefPoint(-5,-1){X}
    \tkzDefPoint(-5/1.3,-1/1.3){Hx}
    \tkzDefPoint(-1,3.5){x}
    
    \tkzDefPoint(-5,2.5){A}
    \tkzDefPoint(-1,.8){B}
    \tkzDefPoint(1.3,-1){C}
    \tkzDefPoint(-3,1){D}
    
    \tkzDefPoint(-2.6,1.1){u1}
    \tkzDefPoint(-2.1,1.9){u2}

   \tkzDrawPolygon[fill=gray!10](A,B,C,D)
    
    \tkzDrawSegments[black,->](O,X)
    \tkzDrawSegments[black,->](O,Y)
    \tkzDrawSegments[black,->](O,Z);
    
    \tkzDrawSegments[blue,->,thick](O,x)
    \tkzLabelSegment[blue,above,pos=1](O,x){$x$}
    
    \tkzDrawSegments[blue,->,thick](O,Hx)
    \tkzLabelSegment[blue,below,pos=.6](O,Hx){$\; Hx = \alpha e_{\cdot 1}=\begin{bmatrix} 
    \alpha\\ 0\\0 \end{bmatrix}$}

    \tkzDrawSegments[red,->](u1,u2);
    \tkzLabelSegment[red,above,pos=1](u1,u2){$u$}
    
    
    \end{tikzpicture}
\end{center}
\captionof{figure}{The mirror $H$ reflects $x$ on $x_1$-axis. The normal vector $u$ is chosen parallel to $x-\alpha e_{\cdot 1}$.}\label{fig:mirror}
\end{center}

Theorem \ref{thm:househ_zero} works with both negative and positive signs. However, to avoid cancelation errors in computing the first component of $u$, one can get rid of subtraction at all by choosing $\mathrm{sign}(x_1)$ in place of $\pm$.
We always form the matrix $H$ via vector
$$
u = x + \mathrm{sign}(x_1)\|x\|_2e_{\cdot 1}, \quad \mbox{with  }\;\mathrm{sign}(0)=+.
$$

Here we discuss how the computational costs of matrix-vector and matrix-matrix multiplications can be reduced when the matrix is a Householder reflection. The standard matrix-vector multiplication $Hx$ requires $2n^2$ flops for a $n\times n$ matrix $H$ and a $n$-vector $x$. The matrix-matrix
multiplication $HA$ for $n\times m$ matrix $A$ costs for $2mn^2$ flops. However, if $H$ is a Householder matrix it is not required forming $H$ explicitly in favour of working with $u$ directly.
In this case, letting $\beta ={2}/{u^Tu}$ one can write
\begin{align*}
Hx & = (I-\beta uu^T)x = x - \beta u(u^Tx)=x - \beta \gamma u, \quad \gamma = u^Tx,\\
HA & = (I-\beta uu^T)A = A - \beta u(u^TA)=A - \beta uw^T, \quad w = A^Tu,\\
AH & = A(I-\beta uu^T) = A - \beta (Au)u^T=A - \beta wu^T, \quad w = Au.
\end{align*}

\begin{workout}
For a Householder matrix $H\in \R^{n\times n}$
verify that the flop-counts for matrix-vector product $Hx$ is about $6n$, and matrix-matrix product $HA$ for $A\in\R^{n\times m}$ (or $AH$ for $A\in\R^{m\times n}$) is about $4mn$. Compare with explicit multiplications.
\end{workout}
\vsp

The following Python code computes a Householder
vector $u$ for a given vector $x$. 
\begin{shaded}
\vspace*{-0.3cm}
\begin{verbatim}
def HouseVec(x):
    # HouseVec(x) computes the Householder vector u such that
    # (I-2uu'/u'u)x = |x|_2e_1
    n = len(x); u = np.zeros(n)
    u[1:] = x[1:]
    s = np.sign(x[1])
    if s == 0: s = 1
    u[0] = x[0] + s*np.linalg.norm(x,2)
    return u
\end{verbatim}
\vspace*{-0.3cm}
\end{shaded}

Multiplication by a Householder transformation is implemented in the following
code:
\begin{shaded}
\vspace*{-0.3cm}
\begin{verbatim}
def HouseProd(u,A):
    # HouseProd(u,A) computes the product of Householder matrix
    # (I-2uu'/u'u) by matrix A
    b = 2/np.dot(u,u); w = np.matmul(np.transpose(A),u)
    HA = A - b*np.outer(u,w)
    return HA
\end{verbatim}
\vspace*{-0.3cm}
\end{shaded}
\vsp
\begin{example}\label{ex:householder_zeros}
The following script transforms the first column of matrix
$$A=\begin{bmatrix}
      2 & 3 & 5 \\
      1 & 2 & -1 \\
      2 & 5 & 3 \\
      1 & -1 & 0
    \end{bmatrix}$$
to a multiple of $e_{\cdot 1}$. We write
\begin{shaded}
\vspace*{-0.3cm}
\begin{verbatim}
A = np.array([[2,3,5],[1,2,-1],[2,5,3],[1,-1,0]])
u = HouseVec(A[:,0])
A = HouseProd(u,A)
print ("Transferred A = \n",  np.round(A,4))
\end{verbatim}
\vspace*{-0.3cm}
\end{shaded}
\noindent
and get the following output (rounded to 4 decimals places):
\begin{shaded}
\vspace*{-0.3cm}
\begin{verbatim}
Transferred A =
 [[-3.1623 -5.3759  -4.7434]
 [  0.      0.3775  -2.8874]
 [  0.      1.755   -0.7749]
 [  0.     -2.6225  -1.8874]]\end{verbatim}
\vspace*{-0.3cm}
\end{shaded}
\noindent

\end{example}
\vsp

\subsection{Plane rotations}
Householder transformations are efficient for dense matrices because of a few flops they need.
In this section we introduce the {\em plane rotations} which are  flexible and can be used efficiently for sparse matrices because they  produce zeros entry by entry. 

The $2\times 2$ skew symmetric matrix
\begin{equation}\label{rotmat22}
  J = \begin{bmatrix}
    c & s \\
    -s & c
  \end{bmatrix}, \quad c^2+s^2=1
\end{equation}
is an orthogonal matrix.
If $c = \cos\theta$ then $Jx$ is a clockwise rotation of $x$ by angle $\theta$. So the matrix $J$
is a rotation matrix in the $(1,2)$-plane.
Sometimes $J$ is called Givens rotation after Wallace
Givens, who used them for eigenvalue computations around 1960. However, they had been used
long before that by Jacobi for the same reason.
Let
$x = (x_1,x_2)^T\neq 0$ and $c = x_1/\|x\|_2$ and $s = x_2/\|x\|_2$ then
$$
Jx = \frac{1}{\|x\|_2}  \begin{bmatrix}
    x_1 & x_2 \\
    -x_2 & x_1
  \end{bmatrix}
   \begin{bmatrix}
    x_1 \\
    x_2
  \end{bmatrix}
  =
   \begin{bmatrix}
    \|x\|_2 \\
    0
  \end{bmatrix}.
$$
In this case, the rotation matrix $J$ rotates $x$ and puts it on $x_1$-axis, i.e. zeros its second component.
By embedding a $2\times2$ rotation in a larger identity matrix, one can manipulate
vectors and matrices of arbitrary dimensions.
\begin{example}
Let $x=(x_1,x_2,x_3,x_4,x_5)^T$ be given such that $\alpha = \sqrt{x_3^2+x_5^2}\neq0$. Let
$c = {x_3}/{\alpha}$ and $s = {x_5}/{\alpha}$. Then we have
$$
\begin{bmatrix}
  1 & 0 & 0 & 0 & 0 \\
  0 & 1 & 0 & 0 & 0 \\
  0 & 0 & c & 0 & s \\
  0 & 0 & 0 & 1 & 0 \\
  0 & 0 &-s & 0 & c
\end{bmatrix}
\begin{bmatrix}
  x_1 \\
  x_2 \\
  x_3 \\
  x_4 \\
  x_5
\end{bmatrix}=
\begin{bmatrix}
  x_1 \\
  x_2 \\
  \alpha \\
  x_4 \\
  0
\end{bmatrix}.
$$
In fact this matrix is a rotation matrix in ($3,5$)-plane which changes $x_5$ to $0$, $x_3$ to $\alpha$ and leaves other components unchanged.
\end{example}
\vsp

Using this idea we can construct a sequence of plane rotations to transfer an arbitrary vector $x\in\R^m$
to a multiple of unit vector $e_{\cdot 1}=[1,0,\ldots,0]^T\in\R^m$.
Let the rotation matrix in the ($j,k$)-plane when applying on vector $x$ with
$c = {x_j}/{\alpha}$ and $s = {x_k}/{\alpha}$ for $\alpha = \sqrt{x_j^2+x_k^2}$
be denoted by $J(j,k)$.
Then it is easy to show that
$$
\underbrace{J(1,2)J(1,3)\cdots J(1,m)}_{G}x = \|x\|_2 e_{\cdot 1}
$$
We note that any rotation matrix $J(j,k)$ with $j<k$ can be used instead of $J(1,k)$ matrices.

\begin{remark}
Possible overflow and underflow in computing $c=x_1/\sqrt{x_1^2+x_2^2}$ and $s=x_2/\sqrt{x_1^2+x_2^2}$ can be avoided by an appropriate scaling. If $|x_1|\geq |x_2|$ we put
$$
t=x_2/x_1, \quad c = 1/\sqrt{1+t^2},\quad s = c\cdot t,
$$
and if $|x_1| < |x_2|$,
$$
t = x_1/x_2, \quad s = 1/\sqrt{1+t^2},\quad c = s\cdot t.
$$
In either case, we can avoid squaring any magnitude larger than $1$.
\end{remark}
\vsp

If an elementary plane rotation $J(j,k)$, with $c= a_{j\ell}/\sqrt{a_{j\ell}^2+a_{k\ell}^2}$ and $s= a_{k\ell}/\sqrt{a_{j\ell}^2+a_{k\ell}^2}$
is applied on a matrix $A$ with entries $a_{j\ell}$ then the $(k,\ell)$ entry of $A$ becomes zero if $k>j$. It is important to note that only two rows $j$ and $k$ of the matrix are changed. This should be taken into account in programming.
Instead of explicitly embedding the $2 \times 2$ rotation matrix into a matrix of larger dimension, which would require unnecessary computations, we can save operations and storage by implementing the rotation more efficiently.
Here are two Python functions that illustrate how to implement the rotation while minimizing computational costs. In the \texttt{GivensProd} function we only operate on two rows of matrix $A$. 

\begin{shaded}
\vspace*{-0.3cm}
\begin{verbatim}
def GivensPar(x,y):
    # GivensPar(x,y) computes Givens parameters to make the second
    # component of [x,y] zero
    if abs(x) > abs(y):
        t = y/x; c = 1/np.sqrt(1+t**2); s = c*t
    else:
        t = x/y; s = 1/np.sqrt(1+t**2); c = s*t
    return c,s
\end{verbatim}
\vspace*{-0.3cm}
\end{shaded}
\begin{shaded}
\vspace*{-0.3cm}
\begin{verbatim}
def GivensProd(c,s,j,k,A):
    # GivensProd(c,s,j,k,A) applies a (j,k)-plane rotation to matrix A
    A[[j,k],:] = np.matmul([[c,s],[-s,c]],A[[j,k],:])
    return A
\end{verbatim}
\vspace*{-0.3cm}
\end{shaded}
\noindent

\begin{example}\label{ex:givens_zeros}
The following script transforms the first column of matrix
$$A=\begin{bmatrix}
      2 & 3 & 5 \\
      1 & 2 & -1 \\
      2 & 5 & 3 \\
      1 & -1 & 0
    \end{bmatrix}$$
to a multiple of $e_{\cdot 1}$.
Consider the matrix $A$ in Example \ref{ex:householder_zeros}.
We use Givens rotations to transfer $A$ into a new matrix that all the components of its first column except the first are annihilated.
\begin{shaded}
\vspace*{-0.3cm}
\begin{verbatim}
A = np.array([[2,3,5],[1,2,-1],[2,5,3],[1,-1,0]])
print("Original A = \n", A)
for k in range(3,0,-1):
    c,s = GivensPar(A[0,0],A[k,0])
    A = GivensProd(c, s, 0, k, A)
print ("Transferred A = \n",  np.round(A,4))
\end{verbatim}
\vspace*{-0.3cm}
\end{shaded}
\noindent
The output with rounding to four decimal places:
\begin{shaded}
\vspace*{-0.3cm}
\begin{verbatim}
Original A =
 [[ 2.  3.  5.]
 [ 1.  2. -1.]
 [ 2.  5.  3.]
 [ 1. -1.  0.]]
Transferred A =
 [[3.1623 5.3759  4.7434]
 [ 0.     0.3162 -2.6352]
 [ 0.     2.2361 -0.7454]
 [ 0.    -2.2361 -2.2361]]
\end{verbatim}
\vspace*{-0.3cm}
\end{shaded}
\noindent
Comparing with the output of the Householder transformation (Example \ref{ex:householder_zeros}), we observe that
$GA$ is not necessarily identical with $HA$. However, in both cases $A$ is transferred to a matrix with zeros under its $a_{11}$ entry.
\end{example}
\vsp 

For a single call the cost of function \verb+GivensPar+ is $6$ flops, and for a matrix $A\in\R^{m\times n}$
the total number of flops in \verb+GivensProd+ is $6n$ (multiplying a $2\times2$ matrix by a $2\times n$ matrix.)
To zero all off-diagonal entries of the first column of a matrix $A\in\R^{m\times n}$, both functions \verb+GivensPar+ and \verb+GivensProd+
should be called in a $m-1$ folds loop. The total cost is thus $(m-1)(6n+6)\approx 6mn$ flops.
\vsp

\subsection{QR factorization}
Matrix decompositions are essential, giving rise to fast and efficient algorithms in matrix computations. Students are typically familiar with LU decompositions (or Gaussian elimination), which factorize a matrix $A$ into a lower triangular matrix $L$ multiplied by an upper triangular matrix $U $, or its other variant $LL^T$ for positive definite matrices. 
Here, we introduce another useful factorization with numerous applications in efficient linear algebra algorithms. 
 
\begin{theorem}
Every matrix $A\in R^{m\times n}$ with $m\geqslant n$ (overdetermined) can be factorized as 
$$
A = QR
$$
where $Q$ is a $m\times m$ {\em orthogonal} matrix and $R$ is a $m\times n$ {\em upper triangular} matrix. 
\end{theorem}
$$
\begin{array}{cccc}
\begin{pmatrix}
  \times & \times & \times  \\
  \times & \times & \times  \\
  \times & \times & \times  \\
  \times & \times & \times  \\
  \times & \times & \times
\end{pmatrix}&=&
\begin{pmatrix}
  \times & \times & \times & \times & \times \\
  \times & \times & \times & \times & \times \\
  \times & \times & \times & \times & \times \\
  \times & \times & \times & \times & \times \\
  \times & \times & \times & \times & \times
\end{pmatrix}&
\begin{pmatrix}
\times & \times  & \times\\
 0 & \times  & \times\\
 0 &  0  & \times \\ \hline
 0 & 0 & 0\\
  0 & 0 & 0
\end{pmatrix}
\\ 
A&& Q&R
\end{array}
$$
\vsp
Here is an illustration for a case with $m=5$ and $n = 3$. 

In the next section we give a ``constructive'' proof for this theorem using Householder reflections. By constructive we mean that the proof also suggests an algorithm to compute the $Q$ and $R$ factors. 

As we observe, the last $m-n$ rows of $R$ are zeros so the last $m-n$ columns of $Q$ have no contribution to the product (but are still important!).  
\vsp
$$
\begin{array}{cccc}
\begin{pmatrix}
  \times & \times & \times  \\
  \times & \times & \times  \\
  \times & \times & \times  \\
  \times & \times & \times  \\
  \times & \times & \times
\end{pmatrix}&=&
\begin{pmatrix}
  \times & \times & \times & \colortxt{$\times$~~~$\times$}  \\
  \times & \times & \times & \colortxt{$\times$~~~$\times$} \\
  \times & \times & \times & \colortxt{$\times$~~~$\times$} \\ 
  \times & \times & \times &\colortxt{$\times$~~~$\times$} \\ 
  \times & \times & \times & \colortxt{$\times$~~~$\times$}
\end{pmatrix}&
\begin{pmatrix}
\times ~~~ \times  ~~~  \times\\
 0 ~~~~ \times  ~~~ \times\\
 0 ~~~~~ 0  ~~~~ \times \\ \hline
 \colortxt{$0$~~~~~$0$~~~~~$0$} \\
 \colortxt{$0$~~~~~$0$~~~~~$0$} 
\end{pmatrix}
\\ \\
A&=& [~~~Q_1~~~~~~~~~Q_2~~]& \begin{bmatrix} R_1\\ 0\end{bmatrix}
\end{array}
$$
This suggests the {\em reduced QR factorization} 
$$
A = Q_1R_1
$$
which is often sufficient for many applications where the $Q_2$ portion is not needed.
\vsp

\subsection{QR factorization using Householder transformations}
Now we show how the idea of introducing zeros in a vector using a Householder matrix can be used for computing a full QR
factorization
$$
A = QR
$$
of a matrix $A\in\R^{m\times n}$, $m\geqslant n$, where $Q\in\R^{m\times m}$ is orthogonal and $R\in\R^{m\times n}$ is upper triangular (entries below the main diagonal are all zero).
This process was first introduced by Householder in 1958. The idea is to reduce $A$ to an upper triangular matrix $R$ by successively
premultiplying $A$ with a series of orthogonal Householder matrices. The products of Householder matrices then constitute the orthogonal matrix $Q$. The process is illustrated for $m=5$ and $n=3$:\\
  Step 1: \vsp
  $$\;H_1A\;\, = H_1\begin{pmatrix}
                                                      \times & \times & \times \\
                                                      \times & \times & \times \\
                                                      \times & \times & \times \\
                                                      \times & \times & \times \\
                                                      \times & \times & \times
                                                    \end{pmatrix}  =
                                                    \begin{pmatrix}
                                                      + & + & + \\
                                                      0 & + & + \\
                                                      0 & + & + \\
                                                      0 & + & + \\
                                                      0 & + & +
                                                    \end{pmatrix} =:A^{(1)}
  $$
  Step 2:
  $$H_2A^{(1)} = H_2\begin{pmatrix}
                                                      + & + & + \\
                                                      0 & + & + \\
                                                      0 & + & + \\
                                                      0 & + & + \\
                                                      0 & + & +
                                                    \end{pmatrix}  =
                                                    \begin{pmatrix}
                                                      + & + & + \\
                                                      0 & * & * \\
                                                      0 & 0 & * \\
                                                      0 & 0 & * \\
                                                      0 & 0 & *
                                                    \end{pmatrix} =:A^{(2)}
  $$
  Step 3:
  $$\qquad H_3A^{(2)} = H_3\begin{pmatrix}
                                                      + & + & + \\
                                                      0 & * & * \\
                                                      0 & 0 & * \\
                                                      0 & 0 & * \\
                                                      0 & 0 & *
                                                    \end{pmatrix}  =
                                                    \begin{pmatrix}
                                                      + & + & + \\
                                                      0 & * & * \\
                                                      0 & 0 & \star \\
                                                      0 & 0 & 0 \\
                                                      0 & 0 & 0
                                                    \end{pmatrix} =:A^{(3)} =: R
  $$
\\
We have $R = A^{(3)}=H_3A^{(2)}=H_3H_2A^{(1)}=H_3H_2H_1A$. If we define $Q^T = H_3H_2H_1$ then $R = Q^TA$ or $A=QR$ where
$Q=H_1H_2H_3$. Remember that $H_k$ are symmetric and orthogonal.

The general case $A\in \R^{m\times n}$ can be treated similarly.
Let $A = (a_{ij})$. In the first step, the Householder matrix $H_1\in \R^{m\times m}$ is build upon the first column of $A$, i.e.,
$$
x = a_{\cdot 1}=\begin{bmatrix}
      a_{11} \\
      a_{2,1}\\
      \vdots \\
      a_{m,1}
    \end{bmatrix}\in \R^m, \quad  u_1 = a_{\cdot 1}-\mathrm{sign}(a_{11})e_{\cdot 1},\quad   H_1 = I-\frac{2}{u_1^Tu_1}u_1u_1^T.
$$
Then $H_1A$ annihilates the components below $a_{11}$. Other entries of $A$ are changed as well. The new matrix is denoted by $$H_1A=:A^{(1)}$$
with entries $a_{ij}^{(1)}$.

In the second step, a Householder matrix $\wt H_2\in \R^{(m-1)\times (m-1)}$ is formed based on the entries from $2$ to $m$ of the second column of $A^{(1)}$, i.e.,
$$
x = a^{(1)}_{2:m,2}=\begin{bmatrix}
      a_{22}^{(1)} \\
      a_{3,2}^{(1)}\\
      \vdots \\
      a_{m,2}^{(1)}
    \end{bmatrix}\in \R^{m-1}, \quad   u_2 = x-\mathrm{sign}(x_1)e_{\cdot 1},\quad   \wt H_2 = I-\frac{2}{ u_2^T u_2} u_2 u_2^T\in \R^{(m-1)\times (m-1)}
$$
Then the Householder matrix $H_2$ is defined by
$$
H_2 = \begin{bmatrix}
        1 & 0 \\
        0 & \wt H_2
      \end{bmatrix}\in\R^{m\times m}.
$$
In the new matrix
$$H_2A^{(1)}=:A^{(2)}$$
the entries below $a_{22}^{(2)}$ become zero. First row and first column of $A^{(2)}$ are identical with that of $A^{(1)}$ due to the special structure of $H_2$ in its first row and column. Specially, the zeros introduced in the pervious step (in the first column) are not destroyed in the current step.

Similarly, in step $k$ a Householder matrix $\wt H_k\in \R^{(m-k+1)\times (m-k+1)}$ is formed based on the entries from $k$ to $m$ of the $k$-th column of $A^{(k-1)}$, i.e.,
$$
x =\begin{bmatrix}
      a_{k,k}^{(k-1)} \\
      a_{k+1,k}^{(k-1)}\\
      \vdots \\
      a_{m,k}^{(k-1)}
    \end{bmatrix}\in \R^{m-k+1}, \quad   u_k = x-\mathrm{sign}(x_1)e_{\cdot 1},\quad   \wt H_k = I-\frac{2}{ u_k^T u_k} u_k u_k^T\in \R^{(m-k+1)\times (m-k+1)}.
$$
Then the Householder matrix $H_k$ is defined by
$$
H_k = \begin{bmatrix}
        I_{k-1} & 0 \\
        0 & \wt H_k
      \end{bmatrix}\in\R^{m\times m},
$$
where $I_{k-1}$ is the identity matrix of size $k-1$.
In the new matrix
$$
H_kA^{(k-1)}=: A^{(k)}
$$
the entries below $a_{kk}^{(k)}$ are all zero. The first block of $H_k$ (the identity block) ensures that the first $k-1$ rows and columns of $A^{(k-1)}$ remains unchanged in $A^{(k)}$. This means that the zeros introduced in all previous steps are not destroyed. For an economic implementation, the matrix-matrix multiplication should only be done on a submatrix of $A^{(k-1)}$ of size $(m-k)\times (n-k)$.

This process is continued until step $n$. In the last step we make zeros below the diagonal in the last column of $A^{(n-1)}$ by multiplying the Householder matrix $H_n$:
$$
H_nA^{(n-1)}=:A^{(n)}.
$$
The resulting matrix $A^{(n)}$ is an upper triangular matrix of size $m\times n$, let us denote it by $R$. We thus have
$$
R = A^{(n)} = H_nA^{(n-1)} = \cdots = \underbrace{H_nH_{n-1}\cdots H_2H_1}_{Q^T}A = Q^TA
$$
where $Q^T$ is an orthogonal matrix because it is the product of $n$ orthogonal Householder matrices. We simply have
$$
A = QR
$$
where $Q = (H_n\cdots H_2H_1)^T=H_1H_2\cdots H_n$.

\begin{remark}[space complexity]
To minimize the storage, $R$ is stored over $A$ in the upper triangular part. Householder matrices are not required to be stored at all. Instead the components $2$ through end of each vector $u_k$ are stored
in the respective positions of $A$ (instead of zeros), and the first components of all $u_k$ vectors are stored in an auxiliary one-dimensional array.
The matrix $Q$, if it is needed explicitly, can be formed in a postprocessing calculation using the stored $u_k$ vectors in a cheap way. See Workout \ref{wo:formQ}.
We note that, in a majority of practical
applications, it is sufficient to have $Q$ in this factored form, and in many applications, $Q$ is
not needed at all.
\end{remark}
\vsp
\begin{remark}
In step $k$ of the algorithm the entries of the submatrix of A containing row $k$
through $m$ and columns $k$ through $n$, denoted by $A(k : m, k : n )$, are updated and stored over the corresponding entries of $A$ via
\begin{equation}\label{QR:updatingH}
\begin{split}
A(k:m,k:n) =&\, (I-\frac{2}{u_k^Tu_k}u_ku_k^T)A(k:m,k:n) \\
=&\, A(k:m,k:n) - \beta u_k u_k^TA(k:m,k:n).
\end{split}
\end{equation}
In a Python code for QR factorization in step $k$ the subroutine \verb+HouseProd+ can be called with input arguments $u_k$ and $A(k:m,k:n)$.
\end{remark}
\vsp 
\begin{workout}\label{wo:formQ}
Show that it requires about $2n^2(m - n/3)$ flops to compute $R$ in the QR factorization
of $A\in\R^{m \times n}$, $m\geqslant n$ using Householder transformations. This cost does not include the explicit construction of $Q$.
Show that it is required about $\frac{4}{3}m^3$ flops to compute $Q$ explicitly.
\end{workout}
\vsp 
The procedure is the same for $m<n$ but is finished after $m-1$ steps. In this case the upper triangular matrix $R$ is of the form
$[R_1~~R_2]$ where $R_1$ is a $m\times m$ upper triangular and $R_2$ is a full matrix of size $m\times (n-m)$. Here is an illustration for
$m=3$ and $n = 5$:

$$
\begin{array}{cccc}
\begin{pmatrix}
  \times & \times & \times & \times & \times \\
  \times & \times & \times & \times & \times \\
  \times & \times & \times & \times & \times
\end{pmatrix}&=&
\begin{pmatrix}
  \times & \times & \times \\
  \times & \times & \times \\
  \times & \times & \times
\end{pmatrix}&
\begin{pmatrix}
  \times & \times & \times & \times & \times \\
  0 & \times & \times & \times & \times \\
  0 & 0 & \times & \times & \times
\end{pmatrix}\\
A&&Q&R
\end{array}
$$
\ \\

The following Python function computes the QR factorization of a $m\times n$ matrix $A$, either $m\geqslant n$ or $m<n$.
The code handles cases where only $R$, $R$ and $u_k$ vectors and both $Q$ and $R$ are demanded. In the second case
Householder vectors $u_k$ are stored in the lower diagonal part of the output matrix and in additional array \verb+u1+.

\begin{shaded}
\vspace*{-0.3cm}
\begin{verbatim}
def qrfac(A, mode = 'Q&R'):
    # qrfac(A, mode) computes the QR factorization of a (m x n) matrix A
    # mode: 'R', 'R&u' and 'Q&R'. The default mode is 'Q&R'
    m,n = np.shape(A)
    s = min(m-1,n)
    if mode == 'R':
        for k in range(s):
            u = HouseVec(A[k:,k])
            A[k:,k:] = HouseProd(u,A[k:,k:])
        return np.triu(A)
    elif mode == 'R&u':
        u1 = np.zeros(s)
        for k in range(s):
            u = HouseVec(A[k:,k])
            A[k:,k:] = HouseProd(u,A[k:,k:])
            A[(k+1):,k] = u[1:]
            u1[k] = u[0];
        return A, u1
    elif mode == 'Q&R':
        A,u = qrfac(A,'R&u')
        Q = np.eye(m)
        for k in range(s):
            Q[k:,:] = HouseProd(np.append(u[k],A[(k+1):,k]),Q[k:,:])
        return np.transpose(Q), np.triu(A)
    else:
        print("Input mode types 'R', 'R&u' or 'Q&R' ")
\end{verbatim}
\vspace*{-0.3cm}
\end{shaded}

\subsection{QR factorization using Givens rotations}
The QR factorization of a matrix $A\in\R^{m\times n}$ can also be simply obtained using Givens rotations in $s = \min\{m-1,n\}$ steps as below:
\begin{itemize}
\item[-]Step 1: form an orthogonal matrix $G_1 = J(1,m)J(1,m-1)\cdots J(1,2)$ such that
$A^{(1)}=G_1A$ has zeros below its $(1,1)$ entry in the first column.
\item[-]Step 2: form an orthogonal matrix $G_2 = J(2,m)J(2,m-1)\cdots J(2,3)$ such that
$A^{(2)}=G_2A^{(1)}$ has zeros below its $(2,2)$ entry in the second column.
\item[] $\vdots$
\item[-]Step $k$: form an orthogonal matrix $G_k = J(k,m)J(k,m-1)\cdots J(k,k+1)$ such that
$A^{(k)}=G_kA^{(k-1)}$ has zeros below its $(k,k)$ entry in the $k$-th column.
\end{itemize}
The final matrix $A^{(s)}:=R$ is upper triangular and the matrix $Q = G_1^TG_2^T\cdots G_s^T$ is orthogonal and
$$
A = QR.
$$
To optimize the storage used, the matrix $R$
is stored over $A$, and $Q$ is formed implicitly out of Givens parameters (for example using Python function \verb+GivensProd+).
\vsp

\begin{labexercise}
Implement a Python function for QR factorization using Givens rotations. Ask the user for different modes. Call your function for some matrices and print the outputs. 
\end{labexercise}
\vsp
\begin{remark}
The QR factorization with Givens rotations requires $3n^2(m-n/3)$ flops. This does not include the computation of $Q$. Compared with the Householder method, this algorithm is about $1.5$ times more expensive.
However when the matrix $A$ is sparse or has a special structure with lots of zeros in its lower triangle, the Givens approach is cheaper.
For example
in several applications (for example in eigenvalue computation) one needs to find the QR factorization of an upper {\bf Hessenberg matrix}. An upper Hessenberg matrix is similar to an upper triangular matrix with additional nonzero elements on the diagonal right below its main diagonal.
Since an upper Hessenberg matrix $A\in \R^{n\times n}$ has at most $(n-1)$ nonzero subdiagonal entries,
the QR factorization of $A$ can be obtained by only $(n- 1)$ Givens rotations. See the following illustration for $n=4$:
$$
\begin{array}{lllllll}
\begin{pmatrix}
  \times & \times & \times & \times \\
  \times & \times & \times & \times \\
  0 & \times & \times & \times \\
  0 & 0 & \times & \times
\end{pmatrix} & \hspace{-.2cm} \xrightarrow{J(1,2)}  & \hspace{-.2cm}
\begin{pmatrix}
  \times & \times & \times & \times \\
  0 & \times & \times & \times \\
  0 & \times & \times & \times \\
  0 & 0 & \times & \times
\end{pmatrix} & \hspace{-.2cm} \xrightarrow{J(2,3)}  & \hspace{-.2cm}
\begin{pmatrix}
  \times & \times & \times & \times \\
  0 & \times & \times & \times \\
  0 & 0 & \times & \times \\
  0 & 0 & \times & \times
\end{pmatrix} & \hspace{-.2cm} \xrightarrow{J(3,4)}  & \hspace{-.2cm}
\begin{pmatrix}
  \times & \times & \times & \times \\
  0 & \times & \times & \times \\
  0 & 0 & \times & \times \\
  0 & 0 & 0 & \times
\end{pmatrix}

\end{array}
$$

\end{remark}
\vsp 
\begin{remark}
It was shown by Wilkinson in 1965 that the computed $\wh Q$ and $\wh R$ with Givens rotations satisfy
$$
\wh R = \wh Q^T(A+E)
$$
where there exists a constant $c$ independent of $m$ and $n$ such that
$$
\|E\|_F\leqslant c\|A\|_F.
$$
This shows that the algorithm is backward stable.
\end{remark}
\vsp 

\subsection{Other algorithms}
The Gram-Schmidt algorithms (classical and modified versions) are alternative algorithms for computing the thin QR factorization of $A$.
See \cite{Trefethen-Bau:1997,Heath:2018}.
\vsp

\subsection{QR factorization for solving the least squares problem}
Solving the linear least squares problems using the normal equations has two
significant drawbacks: (1) Forming $A^TA$ can lead to loss of information,
(2) The condition number $A^TA$ is the square of that of $A$:
$$
\cond_2(A^TA) = [\cond_2(A)]^2.
$$
We illustrate these points in a couple of examples.
\begin{example}\label{ex:ATAdang}
Let
$$
A = \begin{bmatrix}
      1 & 1  \\
      \epsilon & 0  \\
      0 & \epsilon
    \end{bmatrix}
    \vsp
$$
where $\epsilon>0$ is a small real number. Clearly $A$ has full rank and
$$A^TA = \begin{bmatrix}
      1+\epsilon^2 & 1  \\
      1 & 1+\epsilon^2  \\
    \end{bmatrix}
$$
In the double precision floating point arithmetic if we let $\epsilon$ to be smaller that $10^{-8}$ then
$\fl(1+\epsilon^2)=1$ and the computed matrix
$$\fl(A^TA) = \begin{bmatrix}
      1 & 1  \\
      1 & 1  \\
    \end{bmatrix}
$$
is indeed singular.
\end{example}
\vsp
\begin{example}
In the matrix $A$ of Example \ref{ex:ATAdang} assume $\epsilon = 10^{-4}$. Then we can show that
$\cond_2(A)=\sqrt{2}\times 10^4$ while $\cond_2(A^TA)=2\times 10^8$.
\end{example}
\vsp 

In view of the potential numerical difficulties with the normal equations approach,
we need an alternative that does not require formation of the normal system.
In this section we will explain the use of QR factorization for this purpose while in the sequel an alternative approach through SVD will
be discussed.

Consider again the least squares problem
$$
\min_{x\in\R^n}\|Ax-b\|_2
$$
with overdetermined matrix $A\in\R^{m\times n}$.
The QR factorization transforms this linear least squares problem into a
triangular least squares problem having the same solution.
First assume that $A$ has full rank, i.e., $\rankk(A)=n$. Let $A=QR$ be the QR factorization of $A$ and partition
$Q=[Q_1\; Q_2]$ where $Q_1$ consists of first $n$ columns of $Q$, and $R=\begin{bmatrix}R_1 \\ 0  \end{bmatrix}$ where
$R_1\in \R^{n\times n}$ is an upper triangular matrix. In the reduced form $A=Q_1R_1$.
Since $A$ has full rank all diagonal entries of $R_1$ are nonzero, so it is nonsingular.
We can write
\begin{align*}
  \|Ax-b\|_2^2 & = \|QR x - b\|_2^2 = \|R x - Q^Tb\|_2^2 =
   \left\| \begin{bmatrix}R_1 \\ 0  \end{bmatrix}x -
   \begin{bmatrix}Q_1^Tb \\ Q_2^Tb  \end{bmatrix}
    \right\|_2^2\\
    & = \| R_1 x - Q_1^Tb\|_2^2 + \|Q_2^Tb\|_2^2.
\end{align*}
The minimum is obtained if the first norm on the right-hand side is vanished, i.e.,
$$R_1x = Q_1^Tb.$$
 Since $R_1$ is upper triangular and
nonsingular, a simple backward substitution with $\mathcal{O}(n^2)$ flops gives the least square solution $x$.
The reminder then is $$r = \|Ax-b\|_2=\|Q_2^Tb\|_2.$$
If the reminder is not important to us, a reduced QR factorization is enough for obtaining the least squares solution.
\vsp 

\begin{example}
Consider again Example \ref{ex:height_hills}. The least squares solution to height measurements can be computed using the QR factorization
as below.
\begin{shaded}
\vspace*{-0.3cm}
\begin{verbatim}
import numpy as np
A = np.array([[1.,0,0],[0,1.,0],[0,0,1.],[-1.,1.,0],[-1.,0,1.],[0,-1.,1.]])
b = np.array([1237,1941,2417,711,1177,475])
Q,R = qrfac(A)
print('Q =', np.round(Q,4), '\n R =', np.round(R,4))
x = BackSub(R[0:3,:],Q[:,0:3].T@b)
print('x =',np.round(x,4))
\end{verbatim}
\vspace*{-0.3cm}
\end{shaded}
Note that we also called the backward substitution algorithm for solving upper triangular systems using the Python function \texttt{BackSub}. The code for this function is left as an exercise for the reader.
The outputs are
\begin{shaded}
\vspace*{-0.3cm}
\begin{verbatim}
Q = [[-0.5774 -0.2041 -0.3536  0.5113  0.4878 -0.0235]
     [ 0.     -0.6124 -0.3536 -0.4878  0.0235  0.5113]
     [ 0.      0.     -0.7071 -0.0235 -0.5113 -0.4878]
     [ 0.5774 -0.4082 -0.      0.6664 -0.1786  0.1551]
     [ 0.5774  0.2041 -0.3536 -0.1551  0.6664 -0.1786]
     [ 0.      0.6124 -0.3536  0.1786 -0.1551  0.6664]]
 R =[[-1.7321  0.5774  0.5774]
     [ 0.     -1.633   0.8165]
     [ 0.      0.     -1.4142]
     [ 0.      0.      0.    ]
     [ 0.      0.      0.    ]
     [ 0.      0.      0.    ]]
x = [1236. 1943. 2416.]
\end{verbatim}
\vspace*{-0.3cm}
\end{shaded}

\end{example}
\vsp

When $A$ is rank-deficient, i.e., $\rankk(A)<n$, the solution of the least squares problem in not unique; the problem has infinite number of solutions. In this case, the QR factorization of $A$
still exists, but the upper triangular factor $R$ is singular.
This situation usually arises from a poorly designed experiment,
insufficient data, or an inadequate model. If one insists on forging ahead as is,
a variation of QR factorization with {\em column pivoting} (next section) can be used to find all the solutions. See section \ref{sect-lssvd} for an alternative approach via SVD.
Dealing with rank deficiency also enables us to handle
underdetermined problems, where $m < n$, since the columns of $A$ are necessarily
linearly dependent in that case.
\vsp 

\subsection{QR factorization with column pivoting}

In computing the QR factorization, in each step one can
interchange the column having the maximum Euclidean norm in the submatrix with the pivot column.
This kind of column pivoting makes column interchanges so
that the zero pivots are moved to the lower right hand corner of $R$.
The resulting
factorization then is suitable for solving rank-deficient least squares problems. Let's describe the algorithm step by step.

In step 1, we compute the $2$-norm of all columns of $A$, and interchange the column having maximum norm with the first column by multiplying $A$ with a permutation matrix $P_1$ from right. Then we apply the first step of basic QR factorization on $AP_1$ to transfer its first column to the form $[\alpha_1,0,\ldots,0]^T$. By either Householder or Givens transformations we have
 \begin{equation}\label{qr:pivot-step1}
 Q_1AP_1 = \left[\begin{array}{c|ccc}
             \alpha_1 & \tilde a_{11} & \cdots & \tilde a_{1n} \\ \hline
             0 &  \tilde a_{22}  & \cdots & \tilde a_{2n}  \\
             \vdots & \vdots &  & \vdots \\
             0 &  \tilde a_{m2}  & \cdots & \tilde a_{mn}
           \end{array}\right]:=A^{(1)}
 \end{equation}

In step 2, we compute the $2$-norm  of all columns of the submatrix obtained by ignoring the first row and
column of $A^{(1)}$ and interchange the column having maximum norm with the second column by multiplying $A^{(1)}$ with a permutation matrix $P_2$ from right.
(Note: when the columns are interchanged the full columns should be swapped, not just the portions that lie in the
submatrix). Then we apply the second step of the basic QR factorization on $A^{(1)}P_2$:
 $$
 Q_2A^{(1)}P_2 = \left[\begin{array}{cc|ccc}
             \alpha_1 & \tilde a_{11} & \tilde a_{12} &\cdots & \tilde a_{1n} \\
             0 &  \alpha_2 & \hat a_{23}  & \cdots & \hat a_{2n}  \\ \hline
             0 &  0 & \hat a_{33}  & \cdots & \hat a_{3n}  \\
             \vdots &\vdots &\vdots &  & \vdots \\
             0 &  0 &\hat a_{m2}  & \cdots & \hat a_{mn}
           \end{array}\right]:=A^{(2)}
 $$

We continue in a similar way to higher steps. If the matrix has full rank $n$, the algorithm terminates after $n$ steps where in the final step we have
$$
R = A^{(n)} = Q_nA^{(n-1)}P_n = Q_nQ_{n-1}A^{(n-2)}P_{n-1}P_n = \cdots = Q_nQ_{n-1}\cdots Q_1AP_1\cdots P_{n-1}P_n=: Q^TAP
$$
or
$AP = QR$
where $Q = Q_1^T\cdots Q_n^T$ and $P = P_1\cdots P_n$.
Here $PA$ is indeed a matrix obtained from $A$ by permuting some of its columns.
$R$ is upper triangular and nonsingular.

However, if $A$ is rank-deficient there will come a step at which we are forced to take $\alpha_k=0$. In this step all of the entries of the remaining submatrix are zero.
Suppose this occurs after $r$ steps have been completed and

$$
Q_rQ_{r-1}\cdots Q_1AP_1\cdots P_{r-1}P_r =
\left[\begin{array}{cccc|ccc}
\alpha_1 & \times   & \cdots  & \times   & \times & \cdots & \times \\
0        & \alpha_2 & \cdots  & \times   & \times & \cdots & \times \\
\vdots   &          & \ddots  & \vdots   & \vdots &        & \vdots  \\
0        &   0      & \cdots  & \alpha_r & \times & \cdots & \times\\ \hline
0        &   0      & \cdots  &    0     &   0    & \cdots & 0 \\
\vdots   & \vdots   &         & \vdots   & \vdots &        & \vdots\\
0        &  0       & \cdots  &    0     &   0    & \cdots & 0 \\
           \end{array}\right]=: \begin{bmatrix}
                                  R_{11} & R_{12} \\
                                  0 & 0
                                \end{bmatrix}=R
$$
where $R_{11}\in \R^{r\times r}$ is upper triangular and nonsingular. Its main diagonal entries $\alpha_1,\alpha_2,\ldots,\alpha_r$ are all nonzero. Again we have
$
AP = QR
$
where $Q = Q_1^T\cdots Q_r^T$ and $P = P_1\cdots P_r$.
Since $\rankk(R)=r$ we have $\rankk(A)=r$. This result is summarized in the following theorem.

\begin{theorem}
  Given a matrix $A\in \R^{m\times n}$ with $m\geqslant n$ and $\rankk(A)= r\leqslant n$ there exist a permutation matrix $P\in \R^{n\times n}$, an orthogonal matrix $Q\in \R^{m\times m}$, and an upper triangular matrix $R\in \R^{m\times n}$ of the form
  $R = \begin{bmatrix}
                                  R_{11} & R_{12} \\
                                  0 & 0
                                \end{bmatrix}$
with $R_{11}\in \R^{r\times r}$ and nonsingular such that
$$
AP = QR.
$$
\end{theorem}
\vsp

Now, we address the question of how this factorization can be used to solve the least squares problem. First, we note that the permutation matrix $P$ is orthogonal as it is obtained by swapping the columns of the identity matrix. It is clear that $P^{-1}=P^T$ if is applied from left on a vector $x$ will interchange the corresponding rows in $x$.
Assume that $\rankk (A)=r<n$ and $AP = QR$ is the QR factorization of $A$ with column pivoting. Partition $Q=[Q_1\; Q_2]$ where $Q_1$ consists of first $r$ columns of $Q$.
Using the change of variables $y=P^Tx$ and letting $y = [\wt y,\; \wh y]$ for $\wt y\in \R^{r}$, we can write
\begin{align*}
  \|Ax-b\|_2^2 & = \|APP^Tx-b\|_2^2 = \|QRy - b\|_2^2 = \|R y - Q^Tb\|_2^2 \\
   & =
   \left\| \begin{bmatrix}R_{11} & R_{12} \\ 0 & 0 \end{bmatrix}\begin{bmatrix}\wt y \\ \wh y  \end{bmatrix} -
   \begin{bmatrix}Q_1^Tb \\ Q_2^Tb  \end{bmatrix}
    \right\|_2^2= \| R_{11} \wt y + R_{12}\wh y - Q_1^Tb\|_2^2 + \|Q_2^Tb\|_2^2. \vsp
\end{align*}
There are
many choices of $y = [\wt y,\; \wh y]$ for which the first term in the right hand side is zero. The second term is independent of $y$ (and thus $x$) and determines the reminder of the least squares solutions.
Recall that $R_{11}$ is nonsingular. For any choice
of $\wh y\in \R^{n-r}$ there exists a unique $\wt y\in \R^{r}$ such that
$$
R_{11}\wt y = Q_1^Tb - R_{12}\wh y.
$$
Since $R_{11}$ is upper triangular, $\wt y$ can be calculated using a backward substitution. Finally, 
$$
x = Py
$$
for $y = [\wt y,\; \wh y]$ is a solution to the least squares problem. Since $\wh y$ is arbitrary, we obtain infinite number of least squares solutions $x$. 
\vsp

\begin{remark}
In practice we often do not know the rank of $A$ in advance. After $r$ steps of the
QR factorization with column pivoting, $A$ will have been transformed to the form
$\begin{bmatrix}
   R_{11} & R_{12} \\
   0 & R_{22}
\end{bmatrix}$
where $R_{11}$ is nonsingular. If $\rankk(A) = r$, then in principle $R_{22} = 0$
and the algorithm terminates. However, in the presence of roundoff errors $R_{22}$
is not exactly zero.
In a practical computation we might assign a {\bf numerical rank} $r$ to $A$ if the norm of the largest column of $R_{11}$ is less than a prescribed tolerance. The tolerance should be a small value depending on and the accuracy of the entries and the scale of the original matrix. It is usually set to be
$\delta = 10^{-t}\| A \|_\infty$ where the entries of $A$ are correct to $t$ decimal digits.
This approach generally works well, but unfortunately it is not
$100\%$ reliable, because there exists some nearly rank-deficient triangular matrix with relatively large on-diagonal entries.
Search for the well-known {\em Kahan's matrix} as an example.
We will address a more reliable approach to detect the rank deficiency using SVD in a forthcoming section.
\end{remark}
\vsp
\begin{remark}
Another issue that is worth to be addressed is that
if the norms of the columns are computed in the straightforward
manner at each step then the total cost is about $mn^2 - n^3/3$
flops. This cost can be reduced substantially for steps $2,3,\ldots, r$ by using
information from previous steps. For example let $\kappa_1,\kappa_2,\ldots,\kappa_n$ denote the squares of the norms of columns of $AP_1$
in the first step (thus $\kappa_1$ the largest value). Computing all $\kappa_j$ values counts $2mn$ flops.
For the second step of factorization recall
\eqref{qr:pivot-step1}. Since $Q_1$ is orthogonal, the Euclidian norm of columns of $A^{(1)}$ are identical with that of $AP_1$. So the squares of norms of the submatrix can be calculated as
$$
\kappa_j^{(1)} = \kappa_j - \tilde a_{1j}^2, \quad j=2,3,\ldots,n
$$
using $2(m-1)$ flops instead of the $2(m - 1)(n -1)$ flops that is required using the straightforward calculation.

One can continue to other steps similarly, but before starting the new step $k$, the corresponding column swapping should be applied on $\kappa^{(k-1)}$ as well. The total cost then is
$$
2mn + 2\sum_{k=2}^{r}(m-k)  = 2mn + 2m(r-1) - r(r+1) + 2 \approx 2m(n+r) -r^2.
$$
For case $r=n$ the total cost is about $4mn - n^2$ which shows a remarkable reduction in the cost of the
straightforward algorithm.
\end{remark}
\vsp

\begin{workout}
Show that after the QR factorization with column pivoting, the main
diagonal entries of $R_{11}$ satisfy $|\alpha_1|\geqslant |\alpha_2|\geqslant \cdots \geqslant |\alpha_r|$.
\end{workout}
 \vsp

\begin{labexercise}
Develop a Python function for QR factorization with column pivoting. Considers all the aspects described above. 
Then call your function to factorize some specific matrices, and verify your output by checking the equality $PA=QR$.  
Finally use your function for solving the least squares problem $Ax\cong b$ for a given matrix $A$ and right-hand side vector $b$. The coefficient matrix is assumed to be of either full rank or rank-deficient. Assume that the
matrix rank is unknown, so compute it numerically by considering the on-diagonals of the $R$ factor.
\end{labexercise} 
\vsp

%% file: projectors.tex
\section{Projectors*}\label{sect:projectors}

In Figure \ref{fig:bestapp} we observed that the vector $y = Ax\in  \range(A)$ closest to $b$ in the $2$-norm is the
{\em orthogonal projection} of $b$ onto $\range(A)$.
This observation leads to an algebraic characterization of least squares solutions via the notion of projectors\footnote{The reader can skip this section, as the following sections are independent of its content.}.
\begin{definition}
A projector is a square matrix $P$ that satisfies
\begin{equation*}\label{def:projector}
  P^2 = P.
\end{equation*}
\end{definition}

\begin{center}
\begin{center}
\begin{tikzpicture}[scale = 0.8]
    \tkzDefPoint(0,0){A}
    \tkzDefPoint(6,0){B}
    \tkzDefPoint(10,2.5){C}
    \tkzDefPoint(4,2.5){D}
    \tkzDefPoint(4.5,5){E}
    \tkzDefPoint(7,1.75){F}
    \tkzDrawPolygon[thick,fill=gray!10](A,B,C,D)
    \node[color=black] at (2,0.5) {$\mathrm{range}(P)$};
    \filldraw  (4.5,5) circle (2pt);
    \filldraw  (7,1.75) circle (2pt);
    \node[color=black] at (4.8,5) {$v$};
    \node[color=black] at (7.5,1.75) {$Pv$};
    \draw[->,dashed,color=red,thick] (4.5,5)--(7,1.75)node[right] {};
    \tkzLabelSegment[above,pos=.4,sloped](E,F){$Pv-v$}


\end{tikzpicture}
\end{center}
\captionof{figure}{A non-orthogonal projector}\label{fig:non-orth-project}
\end{center}

This definition includes both
orthogonal and nonorthogonal projectors. In Figure \ref{fig:non-orth-project},
the vector $v\in\R^3$ is projected (nonorthogonal) to two-dimensional subspace $\range(P)$.
In this figure, $Pv$ is
the shadow projected by vector $v$ if one shines a light from the north-west direction onto the subspace $\range(P)$.
It is clear that if $v\in \range(P)$ then $Pv=v$ (i.e. $v$ lies exactly on its own shadow).
In fact, if $v\in \range(P)$ then there exists a vector $x$ such that $v = Px$ and
$Pv = P^2x = Px = v$.
We also observe from the figure that if $Pv\neq v$ then the direction in which the light shines is $Pv-v$.
Applying $P$ on the light direction we obtain
$$
P(Pv-v)=P^2v-Pv = Pv-Pv = 0
$$
which means $Pv-v\in\nul(P)$. The direction of the light depends on $v$ but it is always described by a vector in $\nul(P)$.

If $P$ is a projector then $(I - P)^2=I-2P+P^2=I-P$ which means that $(I-P)$ is also a projector. It is called the {\em complementary projector} to $P$. The matrix $I-P$ projects onto $\range(I-P)$ or equivalently onto $\nul(P)$ because:
\begin{lemma}\label{lem:project1}
If $P$ is a projector then
$\range(I-P) = \nul(P)$ and  $\nul(I-P) = \range(P)$.
\end{lemma}
\proof
If $v\in\nul(P)$ then $Pv=0$ so $(I-P)v=v$ which means $v\in \range(I-P)$.
Conversely, for any $v$, we have
$(I -P)v = v - Pv \in \nul(P)$.
By writing $P = I - (I - P)$ we derive the complementary fact $\nul(I-P) = \range(P)$.
$\qed$

\begin{lemma}
If $P$ is a projector then $\range(P) \cap \nul(P) = \{0\}$.
\end{lemma}
\proof
 Any vector $v$ in both sets $\nul(I-P)$ and $\nul(P)$
satisfies $v - Pv =0$ and $Pv=0$ which gives $v=0$. So, $\nul(I - P) \cap \nul(P) = \{0\}$. Then the result follows by Lemma \ref{lem:project1}.
$\qed$
\vsp

We conclude that any projector $P\in\R^{m\times m}$ separates $\R^m$ into two spaces $\range(P)$ and $\nul(P)$ in the sense that
any vector $v\in \R^m$ can be decomposed to $v = x+y$ where $x\in\range(P)$ and $y\in\nul(P)$. Indeed $x=Pv$ and $y=(I-P)v$ because $v=Pv+(I-P)v$.
In this scenario we may write
$$
\R^m = \range(P)+\nul(P).
$$
Conversely, if $S_1$ and $S_2$ are two subspaces of $\R^m$ such that $S_1\cap S_2=\{0\}$ and
$\R^m=S_1+S_2$ then there is a projector $P$ such that $\range(P) =
S_1$ and $\nul(P) = S_2$. We say that $P$ is the projector onto $S_1$ along $S_2$.

\begin{definition}
A projector $P$ is called an {\bf orthogonal projector} if it is symmetric, i.e. $P=P^T$.
\end{definition}
In Figure \ref{fig:orth-project} an orthogonal projection is illustrated.

\begin{center}
\begin{center}
\begin{tikzpicture}[scale = 0.8]
    \tkzDefPoint(0,0){A}
    \tkzDefPoint(6,0){B}
    \tkzDefPoint(10,2.5){C}
    \tkzDefPoint(4,2.5){D}
    \tkzDefPoint(5,5){E}
    \tkzDefPoint(5,1){F}
    \tkzDrawPolygon[thick,fill=gray!10](A,B,C,D)
    \node[color=black] at (2,0.5) {$\mathrm{range}(P)$};
    \filldraw  (5,5) circle (2pt);
    \filldraw  (5,1) circle (2pt);
    \node[color=black] at (5.3,5) {$v$};
    \node[color=black] at (5.5,1) {$Pv$};
    \draw[dashed,color=red,thick] (5,5)--(5,1)node[right] {};
    \tkzLabelSegment[above,pos=.4,sloped](E,F){$Pv-v$}

    \draw[color=red,thick] (4.75,.98)--(4.75,1.25) node[black, right] {};
    \draw[color=red,thick] (4.75,1.25)--(5,1.30) node[black, right] {};
\end{tikzpicture}
\end{center}
\captionof{figure}{An orthogonal projection}\label{fig:orth-project}
\end{center}

 For an orthogonal projector $P$
 the subspaces $\range(P)$ and $\nul(P)$ are orthogonal because
 the inner product between a vector $Px\in \range(P)$ and a vector $(I-P)y\in \range(I-P)=\nul(P)$ is zero:
 $$
 (Px)^T(I-P)y = x^TP^T(I-P)y = x^TP(I-P)y = x^T(P-P^2)y = 0.
 $$
 We use the notations $P_{\perp}:=(I-P)$ and $\range(P)^{\perp}:=\nul(P)$. Any vector
 $v\in\R^m$ then can be expressed as the sum
 $$
 v = (P+(I-P))v = Pv + P_{\perp}v
 $$
 of mutually orthogonal vectors one in $\range(P)$ and the other in $\range(P)^{\perp}$. We also have the
 Pythagorean relation (prove it!)
 $$
 \|v\|_2^2 = \|Pv\|_2^2  +\|P_\perp v\|_2^2.
 $$

 This concept can be applied to find the solution of the least square problem \eqref{leastsquaredef}. If $P$ is an orthogonal projector to
 $\range(A)$ (find $P$ such that $\range(P)=\range(A)$) then we have
 \begin{equation}\label{LS_project_sol}
 \begin{split}
 \|b-Ax\|_2^2 &= \|P(b-Ax) + P_{\perp}(b-Ax)\|_2^2 \\
 & = \|P(b-Ax)\|_2^2 + \|P_\perp(b-Ax)\|_2^2\\
 & =\|Pb-Ax\|_2^2 + \|P_\perp b)\|_2^2.
 \end{split}
 \end{equation}
  The last equality satisfies because $PA=A$ and $P_\perp A=0$. Since the second norm on the right does not depend on
  $x$, the residual is minimized by minimizing the first norm. But the first norm is minimized by the ideal solution $x$ satisfying
  the overdetermined, but consistent system
  \begin{equation}\label{LS-AxPb}
  Ax = Pb.
  \end{equation}
  If fact, such $x$ exists because $Pb\in \range(P)=\range(A)$. If we multiply both sides by $A^T$
  we have
  $A^TAx = A^TPb = A^TP^Tb = (PA)^Tb = A^Tb$ giving
  $$
  A^TAx = A^Tb,
  $$
  which is the normal equation we already derived. The norm $\|P_\perp b)\|_2$ in \eqref{LS_project_sol} is the norm of residual of the least squares solution.

  How can we construct the projection $P$ explicitly? 
  If $A$ is a full-rank matrix then $A^TA$ is nonsingular and
  \begin{equation*}
    P:=A(A^TA)^{-1}A^T
  \end{equation*}
  is an orthogonal projector to $\range(A)$.
  Because it is symmetric and idempotent and $\range(P)=\range(A)$ (prove!). This means that the vector $y\in \range(A)$ closest to $b$ is
  $$
  \tilde b = Pb = A(A^TA)^{-1}A^T=Ax
  $$
  where $x$ is the solution of least squares problem given by the normal equation. We can write $b$ as a sum
  $$
  b = Pb + P_\perp b = Ax + (b-Ax)  = \tilde b +r
  $$
of two mutual orthogonal vectors $\tilde b \in\range(A)$ and $r\in \range(A)^{\perp}$.
\vsp 

\begin{example}
Consider Examples \ref{ex:height_hills} and \ref{ex:height_hills2}. For solution
$x= [1236, 1943, 2416]^T$, the reminder is
$$
r = b-Ax = [1,-2,1,4,-3,2]^T
\vsp
$$
which is orthogonal to all columns of $A$, i.e., $A^Tr = 0$. The orthogonal projector on to $\range(A)$ is
$$
P = A(A^TA)^{-1}A^T = \frac{1}{4}\begin{bmatrix}
2  & 1 &1 &-1&-1& 0 \\
1  & 2 &1 &1 &0 & -1\\
1  & 1 &2 &0 &1 & 1 \\
-1 & 1 &0 &2 &1 & -1\\
-1 & 0 &1 &1 &2 & 1 \\
0  &-1 &1 &-1&1 & 2
\end{bmatrix}\vsp
$$
and the orthogonal projector on to $\range(A)^{\perp}$ is
$$
P_{\perp}=I-P = \frac{1}{4}\begin{bmatrix}
2  & -1  &-1& 1 &1& 0 \\
-1  & 2  &-1 &-1 &0& 1\\
-1  & -1 & 2 &0 &-1& -1\\
1 & -1 &0 &2 &-1& 1\\
1 & 0 &-1 &-1 &2 &-1\\
0  &1& -1 &1& -1& 2
\end{bmatrix},
$$
so that $b = Pb + P_\perp b = \tilde b + r$.
\end{example}
\vsp

An alternative way to define $P$ is to let $Q\in\R^{m\times n}$ be a matrix whose columns form an
orthonormal basis (i.e., $Q^TQ = I$) for $\range(A)$. Then
$$P := QQ^T$$
is symmetric
and idempotent, so it is an orthogonal projector onto $\range(Q) = \range(A)$.
Then from \eqref{LS-AxPb} we have
$Ax = QQ^Tb$. Multiplying both sides by $Q^T$ gives the square system
$
Q^TAx = Q^T b.
$
We will see later how to compute the matrix $Q$ in such a way that this system
is upper triangular and therefore easy to solve.
\vsp

%% file: lec2_part2.tex
\section{Singular Value Decomposition (SVD)}\label{sect:svd}
Let us continue with one of the most important and practical matrix decompositions in numerical linear algebra.
Our geometric presentation here is motivated by \cite{Trefethen-Bau:1997}.
\subsection{Geometric interpretation}
The SVD of a matrix $A$ can be described by the following geometric fact:
\begin{quote}
The image of the unit sphere $S = \{x\in\R^n: \|x\|_2=1\}$ in $\R^n$ under any matrix $A\in \R^{m \times n}$ is the hyperellipsoid
$E = \{Ax: \|x\|_2=1\}$ in
 $\R^m$.
\end{quote}
Hyperellipsoid is just the $m$-dimensional generalization of an ellipse, i.e.,
the surface obtained by stretching the unit sphere in $\R^m$ by some factors
$\sigma_1,\ldots, \sigma_m$ (possibly zero) in some orthonormal directions $u_1,\ldots,u_m\in \R^m$.
The vectors $\{\sigma_ku_k\}$ are the {\em principal semiaxes} of the hyperellipse, with lengths
$\sigma_1,\ldots,\sigma_m$. If A has rank $r$, exactly $r$ of the lengths $u_k$ will turn out to be
nonzero, and in particular, if $m \geqslant n$, at most $n$ of them will be nonzero.

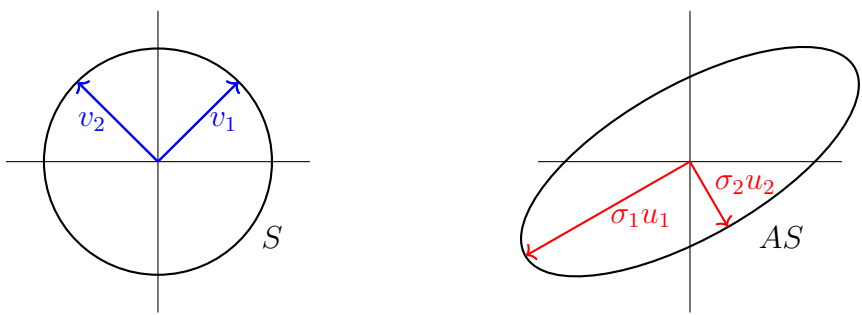
\begin{figure}[!th]
\begin{center}
\begin{tikzpicture}
\draw[black, thick] (-7,0) circle (1.5);
\draw[rotate=30,thick] (0,0) ellipse (2.5 and 1);
\draw (-7,-2)--(-7,2) node[right] {};
\draw (-9,0)--(-5,0) node[right] {};
\draw (0,-2)--(0,2) node[right] {};
\draw (-2,0)--(2,0) node[right] {};
\draw[->,blue,thick] (-7,0)--(-7+1.06,1.06) node[right] {};
\draw[->,blue,thick] (-7,0)--(-7-1.06,1.06) node[right] {};
\draw[->,blue,thick] (-7,0)--(-7+1.06,1.06) node[midway,right] {$v_1$};
\draw[->,blue,thick] (-7,0)--(-7-1.06,1.06) node[midway,left] {$v_2$};
\draw[->,red,thick] (0,0)--(-2.165,-1.25) node[pos=.6,right] {$\hspace{.1cm}\sigma_1u_1$};
\draw[->,red,thick] (0,0)--(0.5,-0.866) node[pos=.35,right] {$\sigma_2u_2$};
\node[color=black] at (-5.5,-1) {$S$};
\node[color=black] at (1.2,-1) {$AS$};
\end{tikzpicture}
\end{center}
\caption{SVD of a $2\times 2$ matrix $A$}
\end{figure}

Let $S$ be the unit sphere in $\R^n$ and $A\in \R^{m\times n}$ with $m\geqslant n$. For simplicity assume that
$A$ has full rank $n$. The image $AS$ is a hyperellipse in $\R^m$. We define $n$ {\bf singular values} of $A$ as the length of $n$ principal semiaxes of $AS$, i.e., $\sigma_1,\sigma_2,\ldots,\sigma_n$. We assume that singular values are numbered in degreasing order,
$$
\sigma_1\geqslant \sigma_2 \geqslant \cdots \geqslant \sigma_n\geqslant 0.
$$
Then we define the {\bf left singular vectors} of $A$ as the direction of principal semiaxes of $AS$. These are orthonormal vectors
$u_1,u_2,\ldots,u_n$ corresponding to singular values $\sigma_1,\ldots,\sigma_n$. This means that  the vector
$\sigma_ku_k$ is the $k$-th largest semiaxis of $AS$.

The {\bf right singular vectors} of $A$ are the orthonormal vectors $v_1,v_2,\ldots,v_n$ that are the preimages of the principal semiaxes of $AS$, numbered so that
\begin{equation}\label{svd:expand}
Av_k = \sigma_k u_k, \quad k=1,2,\ldots,n.
\end{equation}
Equations \eqref{svd:expand} in a compact form can be written as $AV =  U_1 \Sigma_1 $ where $V=[v_1\; v_2\cdots v_n]\in\R^{n\times n}$, $ U_1=[u_1\; u_2\cdots u_n]\in\R^{m\times n}$ and $\Sigma_1 = \mathrm{diag}(\sigma_1,\ldots,\sigma_n)\in \R^{n\times n}$. Since $V$ is an orthogonal matrix, we may write
\begin{equation}\label{svd:reduce}
A =  U_1  \Sigma_1 V^T
\end{equation}
which is called the {\em reduced SVD} of $A$. Here is an illustration for case $m=5$ and $n=3$:

$$
\begin{array}{ccccc}
\begin{pmatrix}
  \times & \times & \times  \\
  \times & \times & \times  \\
  \times & \times & \times  \\
  \times & \times & \times  \\
  \times & \times & \times
\end{pmatrix}&=&
\begin{pmatrix}
  \times & \times & \times  \\
  \times & \times & \times  \\
  \times & \times & \times  \\
  \times & \times & \times  \\
  \times & \times & \times
\end{pmatrix}&
\begin{pmatrix}
\sigma_1 & 0     & 0\\
0      &\sigma_2 & 0\\
0      & 0     & \sigma_3
\end{pmatrix}&
\begin{pmatrix}
\times & \times  & \times\\
\times & \times  & \times\\
\times & \times  & \times
\end{pmatrix}
\\
A&& U_1&\Sigma_1 &V^T
\end{array}
$$
\newline The term {\em reduced} and tilde symbols on $\Sigma$ and $U$ are used to distinguish the factorization \eqref{svd:reduce} from the more standard {\em full SVD}. Since the column of $ U_1$ are $n$ orthonormal vectors in $\R^m$ and $m\geqslant n$ then (unless when $n=m$) they do not form a basis for $\R^m$. We can adjoin an additional $m-n$ columns to $ U_1$ to extend it to an orthonormal matrix $U\in \R^{m\times m}$. For the product to remain unaltered the last $m-n$ columns of $U$ should multiply by zero. So we extend $ \Sigma_1$ by $m-n$ rows of zeros to get the $m\times n$ matrix $\Sigma$. Now we can write the full SVD as
\begin{equation}\label{svd:full}
  A = U\Sigma V^T
\end{equation}
where $U\in\R^{m\times m}$ and $V\in\R^{n\times n}$ are orthonormal matrices and $\Sigma\in\R^{m\times n}$ is a diagonal matrix that caries the nonnegative singular values $\sigma_1,\ldots,\sigma_n$ on its diagonal. The illustration for $m=5$ and $n=3$ is shown below.

$$
\begin{array}{ccccc}
\begin{pmatrix}
  \times & \times & \times  \\
  \times & \times & \times  \\
  \times & \times & \times  \\
  \times & \times & \times  \\
  \times & \times & \times
\end{pmatrix}&=&
\begin{pmatrix}
  \times & \times & \times & \times& \times \\
  \times & \times & \times & \times& \times \\
  \times & \times & \times & \times& \times \\
  \times & \times & \times & \times& \times \\
  \times & \times & \times & \times& \times
\end{pmatrix}&
\begin{pmatrix}
\sigma_1 & 0     & 0\\
0      &\sigma_2 & 0\\
0      & 0     & \sigma_3 \\
0 & 0&0\\
0& 0&0
\end{pmatrix}&
\begin{pmatrix}
\times & \times  & \times\\
\times & \times  & \times\\
\times & \times  & \times
\end{pmatrix}
\\
A&& U& \Sigma &V^T
\end{array}
$$
\newline
If $A$ is rank deficient, $\rankk(A)=r<n$ say, then the factorization \eqref{svd:full} is still valid. The only difference is that now
$\sigma_{r+1} =\cdots = \sigma_n =0$ and $m-r$ columns of $U$ are `silent'.
This means that $r$ singular values and $r$ left singular vectors of $A$ are determined by the geometry of the hyperellipse.
The last $r$ columns of $V$ have no effect in factorization as well.

Let us bring the SVD in a formal definition and prove it mathematically. We prove the existence and uniqueness (under some conditions) of SVD for any complex matrix $A$.
\begin{theorem}\label{thm:svd}
Let $m$ and $n$ be two arbitrary positive integers.
Every matrix $A\in\C^{m\times n}$ has a full SVD of the form
$$
A = U\Sigma V^*
$$
where $U\in\C^{m\times m}$ and $V\in\C^{n\times n}$ are unitary matrices and $\Sigma\in\R^{m\times n}$ is a real diagonal matrix that carries
the singular values $\sigma_1\geqslant\sigma_2\geqslant \cdots \geqslant \sigma_n\geqslant 0$ on its diagonal.
Furthermore, the singular values are uniquely determined and if $A$ is square the left and right singular vectors are uniquely determined up to complex signs.
 \end{theorem}
\proof \cite{Trefethen-Bau:1997}
For proof of existence, let $\sigma_1 := \|A\|_2$. There exists vectors $v_1\in\C^n$ with $\|v_1\|_2=1$ such that
$\|Av_1\|_2 = \sigma_1$. Let $u_1 := Av_1\in \C^m$. Consider any extension of $v_1$ to an orthonormal basis $\{v_1,v_2,\ldots,v_n\}$
for $\C^n$ and any extension of $u_1$ to an orthonormal basis $\{u_1,u_2,\ldots,u_m\}$
for $\C^m$. Suppose that $V_1$ and $U_1$ denote the unitary matrices with columns $v_j$ and $u_j$. We can write
$$
U_1^*AV_1 = \begin{bmatrix}
              u_1^* \\
              \wt U_1^*
            \end{bmatrix} A \begin{bmatrix}
                              v_1 & \wt V_1
                            \end{bmatrix}
            \begin{bmatrix}
              \sigma_1 & u_1^*A \wt V_1 \\
              0 & \wt U_1^* A\wt V_1
            \end{bmatrix}
=:
            \begin{bmatrix}
              \sigma_1 & w^* \\
              0 & B
            \end{bmatrix}=:S,
$$
where $0$ is column vector of dimension $m-1$, $w\in \C^{n-1}$ and $B\in \C^{(m-1)\times(n-1)}$. Now we show $w$ is indeed zero. If fact we have
$$
\left\|\begin{bmatrix}
              \sigma_1 & w^* \\
              0 & B
            \end{bmatrix}
             \begin{bmatrix}
              \sigma_1 \\
              w
            \end{bmatrix}
              \right\|_2=
\sigma_1^2+w^*w \geqslant \sigma_1^2+w^*w =(\sigma_1^2+w^*w)^{1/2}\left\|
             \begin{bmatrix}
              \sigma_1 \\
              w
            \end{bmatrix}
 \right\|_2
$$
which shows that $\|S\|_2\geqslant (\sigma_1^2+w^*w)^{1/2}$. Since $S$ and $A$ are unitarily equivalent we know that
$\|S\|_2=\|A\|_2 = \sigma_1$. This shows that $w=0$, and
$$
U_1^*AV_1 =
            \begin{bmatrix}
              \sigma_1 & 0 \\
              0 & B
            \end{bmatrix}.
$$
The proof is now completed by induction. The proof for case $m=n=1$ is obvious. Let $B = U_2\Sigma_2 V_2^*$ is the full SVD of $B$. We can write
$$
U_1^*AV_1 =
            \begin{bmatrix}
              \sigma_1 & 0 \\
              0 & U_2\Sigma_2 V_2^*
            \end{bmatrix}=
            \begin{bmatrix}
              1 & 0 \\
              0 & U_2
            \end{bmatrix}
            \begin{bmatrix}
              \sigma_1 & 0 \\
              0 & \Sigma_2
            \end{bmatrix}
            \begin{bmatrix}
              1 & 0 \\
              0 & V_2^*
            \end{bmatrix}
$$
which gives the full SVD of $A$ via
$$
A  = \underbrace{U_1 \begin{bmatrix}
              1 & 0 \\
              0 & U_2
            \end{bmatrix}}_{U}
\underbrace{\begin{bmatrix}
              \sigma_1 & 0 \\
              0 & \Sigma_2
            \end{bmatrix}}_{\Sigma}
\underbrace{\begin{bmatrix}
              1 & 0 \\
              0 & V_2^*
            \end{bmatrix}V_1^* }_{V^*},
$$
which completes the proof of existence. On the other hand we have
$$
A^*AV = V\Sigma^T\Sigma ,\quad \mbox{or}\quad A^*A v_k = \sigma_k^2 v_k, \quad k= 1,2,\ldots,n
$$
which shows $\sigma_k^2$ are eigenvalues of Hermitian matrix $A^*A$. This proves the uniqueness of singular values.
$\qed$

For the remainder of this lecture we will assume, without loss of generality,
that $m \geqslant n$, because if $m < n$ we consider the SVD of $A^T$, and if the SVD of $A^T$
is $U\Sigma V^T$, then the SVD of $A$ is $V\Sigma^TU^T$. Besides, we will assume $A$ is a real matrix although all
results can be simply proved for the a complex SVD.
\vsp 
\subsection{Properties of SVD}
Several matrix properties including rank, norms and condition
number can be extracted form SVD. In additions, SVD provides orthonormal
bases for $\range(A)$ and $\nul(A)$ and orthogonal projections onto  $\range(A)$ and $\nul(A)$.

\begin{theorem}\label{thm:svdper}
Let $\sigma_1\geqslant \sigma_2\geqslant \cdots \geqslant \sigma_n\geqslant 0$ be the singular values of $A\in \R^{m\times n}$ with
$m\geqslant n$.
\begin{enumerate}
  \item $\|A\|_2 = \sigma_1$,
  \item $\|A\|_F = \sqrt{\sigma_1^2+\cdots + \sigma_n^2}$,
  \item $\|A^{-1}\|_2=\frac{1}{\sigma_n}$ when $m=n$ and $A$ is nonsingular,
  \item $\mathrm{cond}_2(A)=\frac{\sigma_1}{\sigma_n}$ if $A$ is nonsingular,
  \item $\rankk(A)=$ number of nonzeros singular values.
\end{enumerate}
\end{theorem}
\proof. The proof is left as an exercise. Use the facts that orthogonal matrices preserve norm 2 and norm Frobenius as well as the rank of matrices. See Workout \ref{wo:svdper}.
\vsp
\begin{workout}\label{wo:svdper}
Prove Theorem \ref{thm:svdper}.
\end{workout}
\vsp
The condition number of a square and nonsingular  matrix $A$ is defined by
$
\cond_2(A)=\|A\|_2\|A^{-1}\|_2
$.
We observed that it can be obtained by dividing the largest singular values by the smallest one.
One can generalize this  notion to rectangular matrices by defining
$$
\cond_2(A) := \frac{\sigma_1}{\sigma_n}, \quad A\in \R^{m\times n}
$$
provided that $A$ has full rank ($\sigma_n\neq 0$). The condition number becomes large for nearly rank-deficient matrices, i.e., matrices with very small $\sigma_n$. When $A$ is rank-deficient, i.e., when $\sigma_n=0$, the condition number can be defined to be infinity. However, in practical application (for example in least squares settings)
the true definition of condition number for a rank-deficient matrix $A\in\R^{m\times n}$, $m\geqslant n$ is
\begin{equation}\label{svd:conddef}
\cond_2(A) := \frac{\sigma_1}{\sigma_r},
\end{equation}
where $\rankk(A)=r$. (i.e., $\sigma_{r+1}=\cdots=\sigma_n=0$). This generalizes the definition of condition number for rank-deficient matrices.

\begin{remark}
In presence of nosy data and roundoff errors it is more practical to talk about {\bf numerical rank} rather than just the rank  of a
matrix. In this cases we should set a criterion to accept a computed singular value to be `zero'. For example we may consider the singular values of order $\texttt{eps}$ (machine epsilon) to be zero. However, the norm of the matrix and the errors in data (entries) should also be taken in to account. Here is a practical criterion \cite{Datta:2010, Golub-VanLoan:2010}:
\begin{quote}
  A computed singular value is accepted to be zero if it is less than or equal to $\delta=10^{-t}\|A\|_\infty$ where the entries of $A$ are correct to $t$ decimal digits.
\end{quote}
Using this criterion, $A$ has numerical rank $r$ if the computed singular values $\wh \sigma_1,\ldots,\wh \sigma_n$ satisfy
$$
\wh\sigma_1\geqslant \wh\sigma_2 \geqslant \cdots \geqslant \wh\sigma_r > \delta \geqslant \wh\sigma_{r+1} \geqslant \cdots \geqslant \wh \sigma_n.
$$
This means that to determine a numerical rank of a matrix one needs to count the large singular values only.
\end{remark}
\vsp
As pointed out, SVD provides orthogonal bases for $\range(A)$ and $\nul(A)$ as well as projections to these subspaces. The following theorem reveals this.
\vsp 
\begin{theorem}
Let $A=U\Sigma V^T$ be the SVD of $A\in\R^{m\times n}$ with $m\geqslant n$, and let $\rankk(A)=r$. Partition
$$
U = [U_1\; U_2],\quad V=[V_1\; V_2]
$$
where $U_1$ and $V_1$ consist of first  $r$  columns of $U$ and $V$, respectively. Then
\begin{enumerate}
  \item The columns of $U_1$ form a basis for $\range(A)$.
  \item The columns of $V_2$ form a basis for $\nul(A)$.
  \item $U_1U_1^T$ is a projection onto $\range(A)$.
  \item $U_2U_2^T$ is a projection onto $\nul(A^T)$.
  \item $V_1V_1^T$ is a projection onto $\range(A^T)$.
  \item $V_2V_2^T$ is a projection onto $\nul(A)$.
\end{enumerate}
\end{theorem}
\proof The proofs are straightforward; however, we refer to \cite{Datta:2010} or similar resources. You may skip items 3--6 if you have decided to skip section \ref{sect:projectors}.
$\qed$
\vsp 

\subsection{Computation of SVD}

SVD has a tight connection to the well-known eigendecomposition of symmetric matrices. As we know, if $B\in \R^{n\times n}$ is a symmetric matrix
then all its eigenvalues $\lambda_k$ are real and the corresponding eigenvectors $v_k$ are orthonormal. From $Bv_k = \lambda_k v_k$ for $k=1,\ldots,n$, we can write
$$
B = VDV^T
$$
where $V = [v_1, v_2,\ldots, v_n]$ is orthogonal and $D=\text{diag}\{\lambda_1,\lambda_2,\ldots,\lambda_n\}$ is diagonal. This decomposition is called 
{\bf eigendecomposition}, and has limited applications in practice as it is only available for symmetric matrices\footnote{If $A$ is not symmetric but has a set of independent eigenvectors then it can be decomposed to $A = VDV^{-1}$ for nonsingular matrix $V$. In this case $A$ is called {\em diagonalizable}.}.    

The connection to SVD is obtained by computing $A^TA$ and $AA^T$ as below:
$$
A^TA = (U\Sigma V^T)^T(U\Sigma V^T) = V\Sigma^T\underbrace{U^T U}_{I}\Sigma V^T = V \underbrace{\Sigma^T\Sigma}_{D} V^T = VDV^T
$$
where $D$ is a $n\times n$ diagonal matrix 
$$
D= \Sigma^T\Sigma = \begin{bmatrix}
\sigma_1^2 &   & 0\\ 
&     \ddots  & \\
0&            & \sigma_n^2
\end{bmatrix}.
$$
Since $V$ is orthogonal, $VDV^T$ is the {eigendecomposition} of symmetrix matrix $A^TA$. 
This means that $\sigma_1^2,\ldots, \sigma_n^2$ are eigenvalues and columns of $V$ are corresponding eigenvectors of $A^TA$. 

The same computation for $AA^T$ shows that
$$
AA^T = U\Sigma\Sigma^T U^T = UDU^T
$$
where $D$ is a $m\times m$ diagonal matrix
$$
D= \Sigma\Sigma^T = \begin{bmatrix}
\sigma_1^2 &   &           &  & &0\\ 
&     \ddots   &            & & &\\
&              & \sigma_n^2 & & &\\
&              &            &0& & \\
&             &            & &\ddots& \\
0&             &            & &      &0
\end{bmatrix}
$$
The columns of $U$ are eigenvectors of $AA^T$ corresponding to the same eigenvalues $\sigma_1^2,\ldots, \sigma_n^2$ plus $m-n$ zero eigenvalues $\sigma_{n+1}^2=\cdots=\sigma_m^2=0$. 
 
We conclude that singular values of $A$ are square roots of eigenvalues of $A^TA$. Columns of $V$ are eigenvectors of $A^TA$, and columns of $U$ are eigenvectors of $AA^T$.
\vsp 
\begin{example}
To compute the SVD of 
$$
A = \begin{bmatrix}
3&2\\ 2&3\\ 2 & -2
\end{bmatrix}
$$
we form 
$$
A^TA = \begin{bmatrix}
17 & 8\\ 8 &17 
\end{bmatrix}, \lambda_1 = 25, \; \lambda_2 = 9
$$
and compute its eigenvalues $\lambda_1 = 25$ and $\lambda_2 = 9$.
These give $\sigma_1 = \sqrt \lambda_1 = \sqrt{25}=5$ and $\sigma_2 = \sqrt \lambda_2 = \sqrt{9}=3$. We can also show that eigenvectors of $A^TA$ are
$$
v_1 = \frac{1}{\sqrt 2}
\begin{bmatrix}
1 \\
1
\end{bmatrix}, \; v_2 = \frac{1}{\sqrt 2}
\begin{bmatrix}
1\\
-1
\end{bmatrix}
$$
which result in
$$
V = \frac{1}{\sqrt 2}
\begin{bmatrix}
1& 1\\
1 & -1
\end{bmatrix}.
$$
On the other side, columns of $U$ are eigenvectors of
$$
AA^T = \begin{bmatrix}
13 & 12 & 2 \\
12 & 13 & -2\\
23 & -2 & 8
\end{bmatrix}.
$$
Without additional computations we know that eigenvalues of this matrix are $\lambda_1=25$, $\lambda_2 = 9$ and $\lambda_3 = 0$ (why?). 
However, the eigenvectors of $AA^T$ are different from those of $A^TA$ and can be computed as 
$$
u_1 = \frac{1}{\sqrt 2} \begin{bmatrix}
1 \\
1 \\
0 
\end{bmatrix}, \; u_2 =
\frac{1}{\sqrt{18}}
\begin{bmatrix}
 1 \\
 -1\\
4
\end{bmatrix}, \; u_3 = \frac{1}{3}
\begin{bmatrix}
2\\-2\\ 1
\end{bmatrix}
$$
which show that
$$
U= \begin{bmatrix}
\frac{1}{\sqrt{2}} & \frac{1}{\sqrt{18}} & \frac{2}{3}\\
\frac{1}{\sqrt{2}} & -\frac{1}{\sqrt{18}} &-\frac{2}{3}\\
0 & \frac{4}{\sqrt{18}}& -\frac{1}{3}
\end{bmatrix}.
$$
\end{example}
\vsp 

However, in practice, a different, faster and computationally more stable algorithm is used to compute the SVD factors. The commonly used algorithm is the {\bf Golub-Kahan-Reinsch}
algorithm which consists of two  phases. In the first phase the matrix $A$ is reduced to a bidiagonal
matrix $B$ using Householder transformations. This step is usually  known as the Golub-Kahan
bidiagonal procedure. Then in the second phase the matrix $B$ is further reduced
to diagonal matrix $\Sigma$ using an iterative method that
successively constructs a sequence of bidiagonal matrices $B_k$ such that each $B_k$ has
possibly smaller off-diagonal entries than the previous one. The second procedure is known as the Golub-Reinsch iterative procedure.
It thus makes sense to call the combined two-stage procedure the Golub-Kahan-Reinsch algorithm.
We refer the reader to chapter 10 of \cite{Datta:2010} for details. This algorithm is proved to be computationally stable.

In Python the algorithm is implemented in \verb+numpy.linalg+ (and also in \verb+scipy.linalg+) module which provides both reduced and full SVD for either real or complex matrices. The default command
\verb+numpy.linalg.svd(A)+
is equivalent with
\begin{shaded}
\vspace*{-0.5cm}
\begin{verbatim}
numpy.linalg.svd(A, full_matrices=True, compute_uv=True, hermitian=False)
\end{verbatim}
\vspace*{-0.5cm}
\end{shaded}
If \verb+full_matrices=False+, the output will be the reduced SVD. Additionally, in the case of
\verb+compute_uv=False+, the unitary matrices $U$ and $V$ are not computed and the output is the vector of singular values only.
Finally, by \verb+hermitian=True+, $A$ is assumed to be Hermitian (symmetric if real-valued), enabling a more efficient method for finding singular values. Here is an example.
\begin{shaded}
\vspace*{-0.3cm}
\begin{verbatim}
from numpy.linalg import svd
from numpy import array
A = array([[1, 3, 2],[4,0,-1],[0.5, 2, 1],[1, 1, 1],[2, 1, -2]])
U,S,Vt = svd(A)
print('Full SVD:','\n U =\n',np.round(U,4),'\n Vt =\n',np.round(Vt,4),
      '\n sigma =\n',np.round(S,4))
U,S,Vt = svd(A,full_matrices=False)
print('\n Reduced SVD:','\n U =\n',np.round(U,4),'\n Vt =\n',np.round(Vt,4),
      '\n sigma =\n',np.round(S,4))
\end{verbatim}
\vspace*{-0.3cm}
\end{shaded}

The output of the above script is (note that Python gives $V^T$ instead of $V$):

\begin{shaded}
\vspace*{-0.3cm}
\begin{verbatim}
Full SVD:
 U =
 [[-0.4575 -0.6636 -0.0238 -0.4791  0.3468]
 [-0.6711  0.4798  0.5011 -0.1833 -0.186 ]
 [-0.2784 -0.4014 -0.2015  0.2431 -0.8134]
 [-0.2626 -0.2225  0.2953  0.8113  0.3689]
 [-0.4403  0.3446 -0.7877  0.1398  0.2176]]
 Vt =
 [[-0.8593 -0.5112  0.0186]
 [ 0.3474 -0.6099 -0.7123]
 [ 0.3755 -0.6056  0.7016]]
 sigma =
 [5.149  4.3804 1.5969]

 Reduced SVD:
 U =
 [[-0.4575 -0.6636 -0.0238]
 [-0.6711  0.4798  0.5011]
 [-0.2784 -0.4014 -0.2015]
 [-0.2626 -0.2225  0.2953]
 [-0.4403  0.3446 -0.7877]]
 Vt =
 [[-0.8593 -0.5112  0.0186]
 [ 0.3474 -0.6099 -0.7123]
 [ 0.3755 -0.6056  0.7016]]
 sigma =
 [5.149  4.3804 1.5969]
 \end{verbatim}
\vspace*{-0.3cm}
\end{shaded}

We note that for the reduced SVD the command \verb+U,S,Vt = svd(A,0)+ can be used instead.
\vsp 

\subsection{Solving the least squares problem using SVD}\label{sect-lssvd}
The SVD provides a particularly
flexible method for solving linear least squares
problems of any shape or rank. Consider again the least squares problem
$$
\min_{x\in\R^n}\|Ax-b\|_2
$$
with overdetermined matrix $A\in\R^{m\times n}$. First, we assume that $A$ is full-rank, i.e., $\rankk(A)=n$. Let $A=U\Sigma V^T$ be the SVD of $A$ and partition
$U=[U_1\; U_2]$ where $U_1$ consists of first $n$ columns of $U$, and $\Sigma=\begin{bmatrix}\Sigma_1 \\ 0  \end{bmatrix}$ where
$\Sigma_1\in \R^{n\times n}$ is a square diagonal matrix with $\sigma_k$ on its diagonal.
Since $A$ is full-rank, all singular values are positive and $\Sigma$ is nonsingular.
Furthermore, we use the change of variables
$y = V^Tx$.
 Then, we can write
\begin{align*}
  \|Ax-b\|_2^2 & = \|U\Sigma V^T x - b\|_2^2 = \|\Sigma V^T x - U^Tb\|_2^2 =
   \left\| \begin{bmatrix}\Sigma_1 \\ 0  \end{bmatrix}y -
   \begin{bmatrix}U_1^Tb \\ U_2^Tb  \end{bmatrix}
    \right\|_2^2\\
    & = \| \Sigma_1 y - U_1^Tb\|_2^2 + \|U_2^Tb\|_2^2.
\end{align*}
Since $V$ is nonsingular, the minimization over $x$ is equivalent to minimization over $y$. The minimum is obtained if the first norm on the right-hand side is vanished, i.e. $y = \Sigma_1^{-1}U_1^Tb$. This gives the least squares solution
\begin{equation}\label{eq:leastsquaressol_1}
x = V\Sigma_1^{-1}U_1^Tb = \sum_{j=1}^{n} \frac{u_j^Tb}{\sigma_j}v_j,
\end{equation}
and the reminder is $r = \|Ax-b\|_2=\|U_2^Tb\|_2$.

Now, let $A$ be a rank-deficient matrix, $\rankk(A)=r<n$, say. Partition
$U=[U_1\; U_2]$ now with $U_1$ consists of first $r$ columns of $U$, and
$$
\Sigma = \left[
\begin{array}{ccc|ccc}
{\sigma_1}         &         &           &     &      & \\
                   & \ddots  &           &     &  0   &\ \\
                   &         &{\sigma_r} &     &      & \\
                   \hline
                   &         &          &\colorbox{yellow}{$0$}    &   & \\
                   &   0     &          &      &\ddots& \\
                   &         &          &      &   & \colorbox{yellow}{$0$} \\                                      
                   \hline
                  &          &          &      &      & \\
                  &   0      &          &      &  0   & \\ 
                  &          &          &      &      & 
\end{array}\right] = \begin{bmatrix} 
\wt \Sigma_1 & 0 \\ 0 & 0 
\end{bmatrix}
$$
where $\wt \Sigma_1\in \R^{r\times r}$ is the upper-left square matrix with positive singular values $\sigma_1,\ldots,\sigma_r$ on its diagonal. Also, assume that 
$$
 y = V^Tx =:\begin{bmatrix}
                                          \wt y \\
                                          \wh y
                                        \end{bmatrix}
$$
\vsp
where
$\wt y\in \R^r$. Similar to the first case, we have
\begin{align*}
  \|Ax-b\|_2^2 & = \|U\Sigma V^T x - b\|_2^2 = \|\Sigma V^T x - U^Tb\|_2^2 =
   \left\| \begin{bmatrix}\wt\Sigma_1 & 0 \\ 0 & 0  \end{bmatrix}\begin{bmatrix}\wt y\\ \wh y \end{bmatrix} -
   \begin{bmatrix}U_1^Tb \\ U_2^Tb  \end{bmatrix}
    \right\|_2^2\\
    & = \| \wt \Sigma_1 \wt y - U_1^Tb\|_2^2 + \|U_2^Tb\|_2^2.
\end{align*}
In this case
the least square solution is not unique. Solutions are obtained by putting $\wt y = \wt\Sigma_1^{-1}U_1^Tb$ and letting $\wh y$ arbitrary, and then
\begin{equation}\label{eq:leastsquaressol_2}
x = Vy = \sum_{j=1}^{r} \frac{u_j^Tb}{\sigma_j}v_j + \sum_{j=r+1}^{n}y_j v_j  .
\end{equation}
A least squares solution $x$ with the minimum norm 2 among others is called the {\bf norm-minimal} solution. 
Since $\|x\|_2^2=\|y\|_2^2 =\|\wt y\|_2^2 + \|\wh y\|_2^2 $, the norm-minimal solution of the rank-deficient least squares problem is obtained by setting $\wh y = 0$, i.e. by dropping the second summation in \eqref{eq:leastsquaressol_2}.

\begin{labexercise}
Write a Python function for solving the least squares problem for either full-rank or rank-deficient case using SVD.
Assume that the rank of the coefficient matrix is unknown in advance, so compute it numerically. 
\end{labexercise}
\vsp 

\subsection{Pseudoinverse}

Assuming $A \in \R^{m \times n} $ has rank $r $ with $r \leqslant n $, the SVD of $A$ is represented as
\begin{equation}\label{svd:rank_deficient}
A = [U_1 \;\, U_2] \begin{bmatrix} \widetilde{\Sigma}_1 & 0 \\ 0 & 0 \end{bmatrix} [V_1 \;\, V_2]^T = U_1 \widetilde{\Sigma}_1 V_1^T
\end{equation}
where $U_1$ and $ V_1 $ consist of the first $ r$ columns of $U$ and $V$, respectively, and $\widetilde{\Sigma}_1 \in \R^{r \times r} $ is a diagonal matrix with positive singular values $\sigma_1, \ldots, \sigma_r$ on its diagonal.

When $m = n$ (i.e., $A$ is square) and $r = n$ (i.e., $ A $ is full-rank), we have the standard representation $A = U \Sigma V^T$ with all matrices square and $\Sigma$ nonsingular. In this case, the inverse of $ A $ is computed as
$$
A^{-1} = (U \Sigma V^T)^{-1} = V^{-T} \Sigma^{-1} U^{-1} = V \Sigma^{-1} U^T.
$$
This suggests an algorithm to compute the inverse of square and nonsingular matrices. However, this algorithm is more costly than the usual algorithm based on Gaussian elimination, although it is more stable for ill-conditioned matrices.

When $A$ is rectangular and even rank-deficient, we can generalize this approach to compute the pseudoinverse of $A$. In this case, we use the factorization \eqref{svd:rank_deficient} and define
\begin{equation*}
A^{+} :=  V_1 \wt\Sigma^{-1}U_1^T
\end{equation*}
which is indeed the {\bf Moore–Penrose inverse} of matrix $A$. As a result, we can see from \eqref{eq:leastsquaressol_1} that the least squares solution of a full-rank system $Ax \cong b$ is given by
$$
x = A^{+} b.
$$
This is the counterpart of the exact solution $ x = A^{-1} b $ for square and nonsingular linear system $Ax = b$. Additionally, from \eqref{eq:leastsquaressol_1}, we observe that for the rank-deficient problem, $x = A^{+} b$ is the norm-minimal solution.
\vsp 

\begin{example}
The SVD factors of a matrix $A$ are given by
$$
U = \begin{bmatrix}
\frac{1}{\sqrt 2} & \frac{-1}{\sqrt 2} & 0 & 0 & 0\\
0 & 0 & 0 &1 & 0\\
0 & 0 & -1 & 0 & 0 \\
\frac{1}{\sqrt 2} & \frac{1}{\sqrt{2}} & 0 & 0 & 0\\
0 & 0 & 0 &0 & 1
\end{bmatrix}, \;
\Sigma = \begin{bmatrix}
2\sqrt 3 & 0 & 0 & 0 \\
0 & 2 & 0 & 0\\
0 & 0 & 0 & 0\\
0 & 0 & 0 & 0\\
0 & 0 & 0 & 0
\end{bmatrix}, \;
V = \begin{bmatrix}
\frac{\sqrt 6}{3} & 0 & 0 & \frac{-1}{\sqrt 3}\\
0 & 0 & 1 & 0\\
\frac{1}{\sqrt 6} & \frac{-1}{\sqrt{2}} & 0 & \frac{1}{\sqrt 3}\\
\frac{1}{\sqrt 6} & \frac{1}{\sqrt{2}} & 0 & \frac{1}{\sqrt 3}
\end{bmatrix}.
$$
The pseudoinverse of $A$ is computed as 
\begin{align*}
A^+ = V_1\Sigma_1^{-1}U_1^T = 
\begin{bmatrix}
\frac{\sqrt 6}{3} & 0 \\
0 & 0 \\
\frac{1}{\sqrt 6} & \frac{-1}{\sqrt{2}} \\
\frac{1}{\sqrt 6} & \frac{1}{\sqrt{2}} 
\end{bmatrix}
\begin{bmatrix}
\frac{1}{2\sqrt 3} & 0 \\
0 & \frac{1}{2} 
\end{bmatrix}
\begin{bmatrix}
\frac{1}{\sqrt{2}} & 0 & 0 & \frac{1}{\sqrt{2}} & 0 \\
\frac{-1}{\sqrt{2}} & 0 & 0 & \frac{1}{\sqrt{2}}& 0
\end{bmatrix} 
= \begin{bmatrix}
\frac{1}{6} & 0 & 0 &\frac{1}{6} & 0\\
0 & 0 & 0 & 0 & 0 \\
1 & 0 & 0 & 0 & 0 \\
0 & 0 & 0 & 1 & 0 \\
\end{bmatrix}.
\end{align*}
In addition, if we aim to
solve the least square problem $\min \|Ax-b\|_2$ with $b = [1,1,1,1,1]^T$ for the norm-minimal solution, we write 
$$
x = A^+b = 
 \begin{bmatrix}
\frac{1}{6} & 0 & 0 &\frac{1}{6} & 0\\
0 & 0 & 0 & 0 & 0 \\
1 & 0 & 0 & 0 & 0 \\
0 & 0 & 0 & 1 & 0 \\
\end{bmatrix}
\begin{bmatrix}
1\\ 1\\ 1\\ 1\\ 1
\end{bmatrix}
= \begin{bmatrix}
\frac{1}{3} \\
0 \\
1 \\
1
\end{bmatrix}.
$$
\end{example}
\vsp 

\subsection{Low-rank approximation}
Low-rank approximation is a fundamental concept in linear algebra and numerical analysis that plays a crucial role in various applications. 
In this approach, we seek to represent a given matrix with a lower-rank approximation that captures its essential structure while reducing its size and complexity. SVD enables us to achieve this kind of approximation. 

We start with the following theorem whihc addresses the question of how far is a rank $r$ matrix from a matrix of rank $k<r$.
This theorem is generally known as
the Eckart-Yaung theorem.
\begin{theorem}\label{thm:svdAk}
Let $A=U\Sigma V^T$ be the SVD of $A$, and let $k\leqslant r = \rankk(A)$. Define $A_k = U\Sigma_kV^T$ where
$\Sigma_k=\diagg\{\sigma_1,\ldots,\sigma_k,0,\ldots,0\}$, where $\sigma_1\geqslant\sigma_2\geqslant \cdots \geqslant \sigma_k>0$. Then the followings hold true:
\begin{enumerate}
  \item $A_k$ has rank $k$,
  \item The distance between $A$ and $A_k$ is $\|A-A_k\|_2=\sigma_{k+1}$.
  \item Out of all rank $k$ matrices, $A_k$ is the closet to $A$, that is
  $$
  \min_{\rankk(B)=k}\|A-B\|_2 = \|A-A_k\|_2.
  $$
  \end{enumerate}
\end{theorem}
\proof
The proof of 1 is obvious because $\rankk(A_k) = \rankk(U\Sigma_k V^T) = \rankk(\Sigma_k)=k$. For the second item we have
$$
\|A-A_k\|_2 = \|U(\Sigma-\Sigma_k)V^T\|_2 = \|\Sigma-\Sigma_k\|_2 = \sigma_{k+1}.
$$
To prove item 3, we assume that $B$ is another $m\times n$ matrix of rank $k$. So the null space of $B$ has dimension $n-k$. Consider the space
$S = \spann\{v_1,\ldots v_{k+1}\}$ where $v_j$ are right singular vectors of $A$. The intersection of $\nul(B)$ and $S$ must be nonempty because both spaces are subspaces of $\R^n$ and the sum of their dimensions is greater than $n$. Let $z$ be a normal vector ($\|z\|_2=1$) lying in this intersection. Since $z\in S$ there exist scalers $c_j$ such that $z = c_1v_1 + \cdots + c_{k+1}v_{k+1}$. Since $\{v_j\}$ are orthonormal
we have $|c_1|^2+\cdots + |c_{k+1}|^2=1$. On the other side, since $z\in \nul(B)$ we have $Bz = 0$. So
$$
(A-B)z = Az = c_1Av_1 + \cdots + c_{k+1}Av_{k+1} = c_1\sigma_1u_1 + \cdots + c_{k+1}\sigma_{k+1}u_{k+1}
$$
and since $\{u_j\}$ are orthonormal we have
$$
\|(A-B)z\|_2^2 = |c_1\sigma_1|^2 + \cdots + |c_{k+1}\sigma_{k+1}|^2\geqslant \sigma_{k+1}^2(|c_1|^2+\cdots + |c_{k+1}|^2) = \sigma_{k+1}^2.
$$
This shows that
$$\|A-B\|_2^2=\max_{\|y\|_2=1} \|(A-B)y\|_2^2 \geqslant \|(A-B)z\|_2^2 \geqslant  \sigma_{k+1}^2 = \|A-A_k\|_2^2
$$
which completes the proof.
$\qed$
\vsp

Another representation of SVD is
\begin{equation*}
A = U\Sigma V^T = \sigma_1 E_1 + \sigma_2 E_2 + \cdots + \sigma_n E_n, \quad E_k = u_kv_k^T.
\end{equation*}
Each $E_k$ is a rank $1$ matrix which shows that SVD writes $A$
as a sum of rank $1$ matrices. Since $\sigma_1\geq \sigma_2\geq \cdots \geq \sigma_n$, matrix $A$ is
expressed as a list of its ``ingredients'', ordered by ``importance''. 
Each elementary matrix $E_k$ can be stored using only $m+n$ storage locations.
Moreover, the matrix-vector product $E_kx$ can be formed using only $2n+m$ flops.
As was shown in Theorem \ref{thm:svdAk},
$$
A_k = U\Sigma_k V^T = \sigma_1 E_1 + \sigma_2 E_2 + \cdots + \sigma_k E_k
\vsp
$$
is $k$-rank approximation of $A$ with $2$-norm error $\sigma_{k+1}$. Such an approximation
is useful in data dimensionally reduction, image processing, information retrieval,
cryptography, and numerous other applications.

%% file: lec2_part3.tex
\vsp
\subsection{Some applications  of SVD}
We have already studied the use of SVD for computing some norms, condition number, orthogonal bases for some subspaces
related to a matrix, and an application for solving the linear least squares problem. In this section
more concrete examples are given from different applications.

\subsubsection*{Image compression}
\begin{figure}[!th]
\centering
\includegraphics[scale=.55]{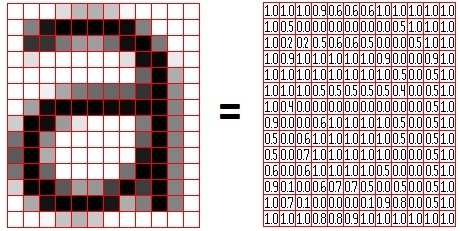}
\caption{Representation of an image with a matrix (form \texttt{pippin.gimp.org/image{\_}processing}).}\label{fig:image_comp}
\end{figure}

An image can be represented by an $m\times n$ matrix $A$ whose $(i,j)$-th
entry corresponds to the brightness of the pixel $(i, j)$.
The storage of this matrix requires $mn$ locations. See Figure \ref{fig:image_comp},

The idea of image compression is to
compress the image represented by a very large matrix to the one which corresponds to a
lower-order approximation of $A$.
SVD provides a simple way if one stores
$$
\sigma_1u_1v_1^T + \cdots + \sigma_ku_kv_k^T =:A_k
\vsp
$$
instead in $(m+n+1)k$ locations (the first $k$ columns of $U$ and $V$ together with the first $k$ singular values).
 This results a considerable savings when $k$ is small.
On the other side $k$ should be large enough to keep the quality of the image still acceptable to the user\footnote{
For image compression, more sophisticated methods like JPG that take human perception into account generally outperform compression using SVD.}.

\begin{figure}[!th]
\centering
\includegraphics[scale=0.6]{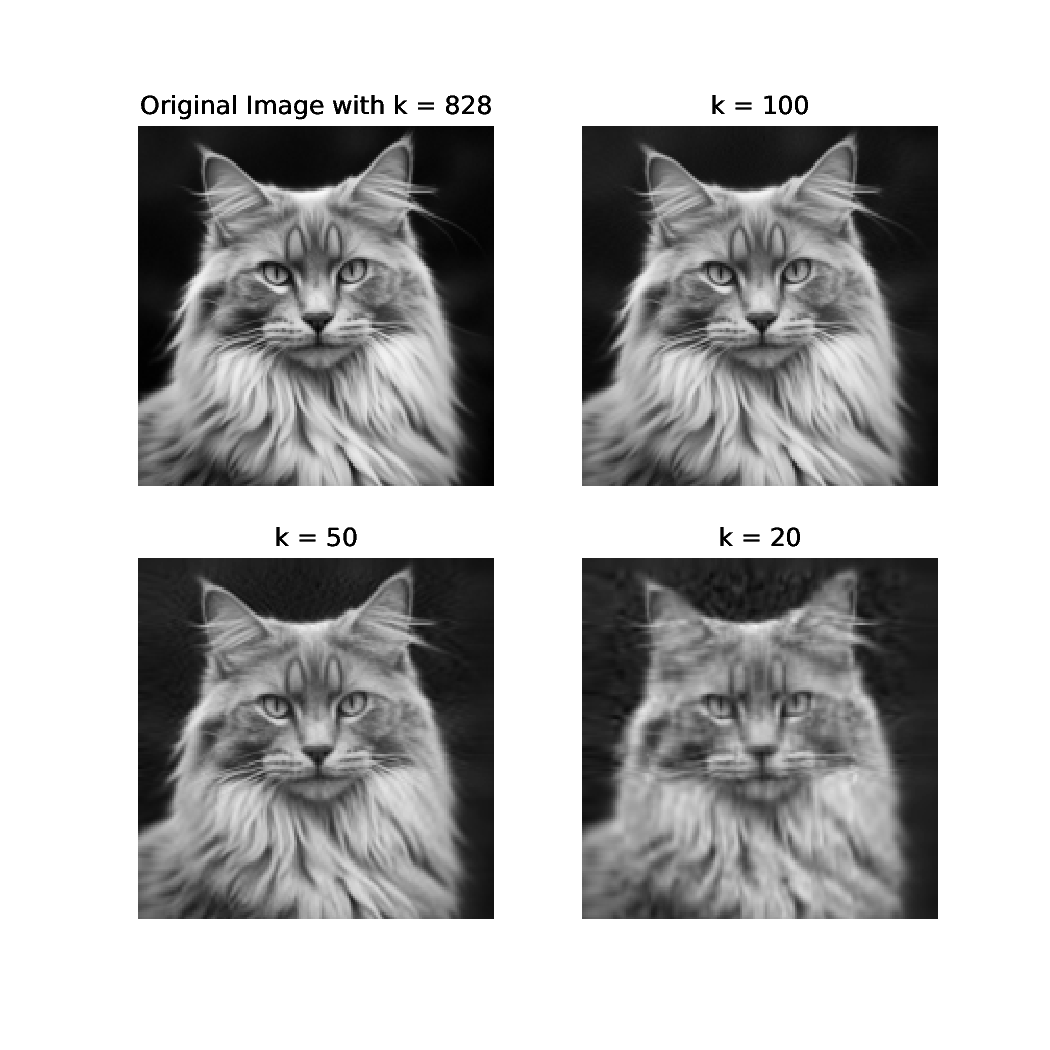}
\caption{Original cat image and compressed images with different values of $k$.}\label{fig:Lena}
\end{figure}

In Figure \ref{fig:Lena} different low-rank approximations of a cat image obtained from the following Python script are shown.

\begin{shaded}
\vspace*{-0.3cm}
\begin{verbatim}
import numpy as np
import matplotlib.pyplot as plt
from PIL import Image

A = Image.open('cat.jpg').convert('L')  # Grayscale read picture
sz = np.shape(A)
U, S, Vt = np.linalg.svd(A)  # Get svd of A
k = [sz[1], 100, 50, 20]
plt.figure(figsize = (7, 7))
for i in range(4):
  Ak = U[:, :k[i]] @ np.diag(S[:k[i]]) @ Vt[:k[i],:]
  if(i == 0):
     plt.subplot(2, 2, i+1), plt.imshow(Ak, cmap ='gray'), plt.axis('off'),
     plt.title("Original Image with k = " + str(k[i]))
  else:
     plt.subplot(2, 2, i+1), plt.imshow(Ak, cmap ='gray'), plt.axis('off'),
     plt.title("k = " + str(k[i]))
plt.savefig("cat1.jpg")
\end{verbatim}
\vspace*{-0.3cm}
\end{shaded}

\subsubsection*{Image restoration}
The aim in image restoration is to restor the original image from a blurry image contaminated by {\em noises}.
It is known that noises correspond to the small singular values.
Thus a simple idea is to compute the SVD of the noisy image and eliminate its last $n-k$ small singular values (smaller than a threshold) and consider the low-rank approximation
$$
A_k = U_k\Sigma_kV_k^T
$$
as a noise-free image.
It is left as an exercise giving an example with a Python script.

Note that, the literature contains a large number of more sophisticated denosing algorithms; some of them are based on SVD.
\vsp

\subsubsection*{Principal component analysis (PCA)}
\begin{figure}[!th]
\centering
\includegraphics[scale=0.55]{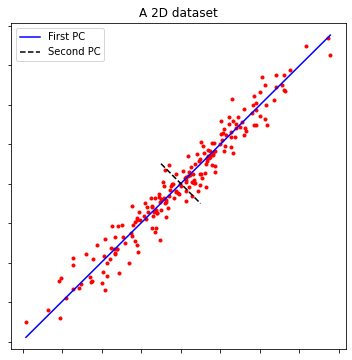}
\caption{Cluster of points in $\R^2$ with principal components (PCs).}\label{fig:pca2D}
\end{figure}

In PCA the objective is to reduce the dimensionality of a dataset in order to use the transferred data in applications that might not work well with high-dimensional data.
The reduction is done by projecting each data point onto only the first few directions ({\em principal components}) to obtain lower dimensional data while preserving as much of the data's variation as possible.
The first principal component is a direction that maximizes the {\em variance} of the projected data.  This direction captures most information of the data.
The second principal component is then calculated with the condition that it is uncorrelated with (i.e., perpendicular to) the first principal component and that it accounts for the next highest variance.
We may proceed through other components, similarly. See Figure \ref{fig:pca2D}.

Let us outline the fundamental principles of PCA within a mathematical framework\footnote{For a history, a review and some recent
developments of PCA see the following sources:
\newline Ian T. Jolliffe, Jorge Cadima, Principal component analysis:
a review and recent developments, Phil. Trans. R. Soc. A 374:
20150202, (2016).\newline
Ian T. Jolliffe, Principal component analysis, 2nd ed. New York, NY, Springer-Verlag, (2002).}.
Assume that we are given a dataset with
observations on $d$ variables (dimensions), for each of $n$ entities or individuals. These data values
define $n$ vectors $x_{\cdot1}, \ldots , x_{\cdot n}$ each in $\R^d$ or, equivalently, a $d \times n$ {\em data matrix}
$$
X=[x_{\cdot1},x_{\cdot2},\ldots,x_{\cdot n}]\in \R^{d\times n}.
$$
We assume that $X$ is row-centred, i.e., the mean value of each row is zero. Otherwise, subtract each row by its mean.
We are looking for a new orthonormal basis $\{v_{\cdot 1},v_{\cdot2},\ldots,v_{\cdot d}\}$, $v_{\cdot j}\in\R^n$, for the column space of $X^T$ (row space of $X$) which
determines the directions of maximum variances in order of significance (a decreasing order). Any element of the new basis is a
linear combination of 
columns of $X^T$. Such linear combinations can be written as
\begin{equation}\label{pca-svdrep}
\sigma_j v_{\cdot j} = \sum_{\ell=1}^{d} u_{\ell j} x_{\ell \cdot} =X^T u_{\cdot j}, \quad j=1,2,\ldots,d,
\end{equation}
where $u_{\cdot j}=[u_{1j},\ldots,u_{dj}]^T$ are assumed to be of the unit norm, i.e. $\|u_{\cdot j}\|_2=1$.
The normalization constants $\sigma_j$ are multiplied from left to make this assumption meaningful.
The variance of any such linear combination
is given by
$$
\varr( X^Tu) = \frac{1}{n-1}(X^Tu)^T(X^Tu) = u^TCu,   
$$
where $C=\frac{1}{n-1}XX^T\in \R^{d\times d}$ is the sample {\em covariance matrix} associated with the dataset.
Hence, identifying the linear combination with maximum variance is
equivalent to obtaining a $d$-dimensional vector $u$ with $u^Tu=1$ which maximizes the quadratic form
$u^TCu$.
By introducing a Lagrange multiplier $\lambda$, the problem is equivalent to maximizing
$$
u^TCu - \lambda(u^Tu-1).
$$
Differentiating with respect to $u$, and equating to the null vector, produces the
equation
\begin{equation*}
  Cu = \lambda u.
\end{equation*}
This means that $u$ is an eigenvector and $\lambda$ its corresponding eigenvalue of the covariance
matrix $C$. Since $C$ is symmetric and semi-positive definite, all eigenvalues are real and nonnegative and eigenvectors form an orthonormal basis for $\R^d$. The eigenvalue $\lambda$ corresponding to eigenvector $u$ of $C$ is indeed the variance of the linear combination defined by $u$, i.e. $X^Tu$, because
$$\varr(X^Tu)=u^TCu=u^T(\lambda u)=\lambda u^Tu = \lambda.$$
Let us denote the pair of eigenvalues and eigenvectors of $C$ by $(\lambda_j, u_{\cdot j})$ for $j=1,2,\ldots, d$ and sort eigenvalues is such way that $\lambda_1\geqslant \lambda_2\geqslant \cdots \geqslant \lambda_d\geqslant 0$.
The maximum variance is attained at the dominant pair $(\lambda_1,u_{\cdot 1})$. Recalling \eqref{pca-svdrep}, this shows that
the maximum variance of the data happens across the direction $u_{\cdot 1}$ with variance
$$
\lambda_1 = \varr(X^Tu_{\cdot1})=\varr(\sigma_1v_{\cdot1}) = \frac{\sigma_1^2}{n-1}v_{\cdot1}^Tv_{\cdot 1} = \frac{\sigma_1^2}{n-1}.
$$
This means that $u_{\cdot1}$ is the first principle component of data $X$ and $\frac{\sigma_1^2}{n-1}$ is the highest variance in the data along the direction $u_1$.

A Lagrange multipliers approach, with the
added restrictions of orthogonality of different vectors $u_{\cdot j}$ can be used to show that
the full set $\{u_{\cdot 1},\ldots,u_{\cdot d}\}$ are the solutions to the problem of obtaining up to $d$ new linear combinations
$X^Tu_{\cdot j}$
which successively maximize the variance, subject to uncorrelatedness
with previous linear combinations. Uncorrelatedness results from the fact that the covariance
between two such linear combinations is zero;
$$
\covv(X^Tu_{\cdot j}, X^Tu_{\cdot k}) = u_{\cdot j}^T C u_{\cdot k} = \lambda_k u_{\cdot j}^Tu_{\cdot k}=0, \quad j\neq k.
$$
As the above formulation reveals, PCA is essentially connected to the SVD. From \eqref{pca-svdrep} we have
\begin{equation}\label{pcr-svd1}
X = U\Sigma V^T
\end{equation}
where $U=[u_{\cdot1}\; u_{\cdot2}\ldots u_{\cdot d}]\in \R^{d\times d}$ and $V=[v_{\cdot1}\; v_{\cdot2}\ldots v_{\cdot d}]\in \R^{d\times n}$ have orthonormal columns and $\Sigma = \diagg(\sigma_1,\sigma_2,\ldots,\sigma_d)$. Consequently, the reduced SVD \eqref{pcr-svd1} gives all principal components as columns of the factor $U$. The variances in each direction are calculated as the squares of the singular values divided by $(n-1)$.

\begin{figure}[!th]
\centering
\includegraphics[scale=0.5]{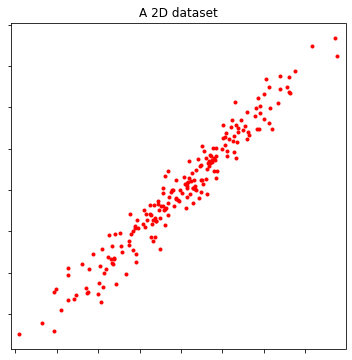} \includegraphics[scale=0.5]{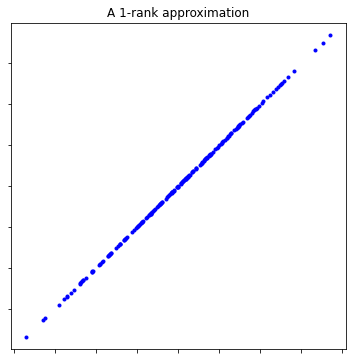}
\caption{A dataset in $\R^2$ (left) and its $1$-rank approximation using SVD.}\label{fig:pca2Dapp}
\end{figure}

\begin{figure}[!th]
\centering
\includegraphics[scale=0.5]{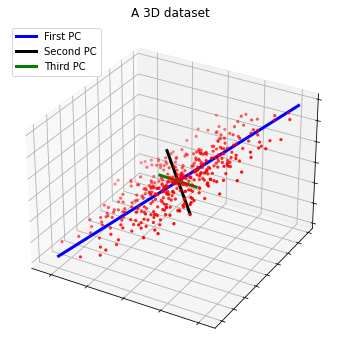} \includegraphics[scale=0.5]{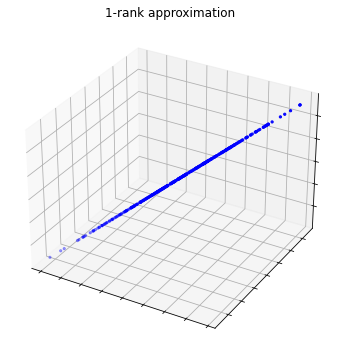}\includegraphics[scale=0.5]{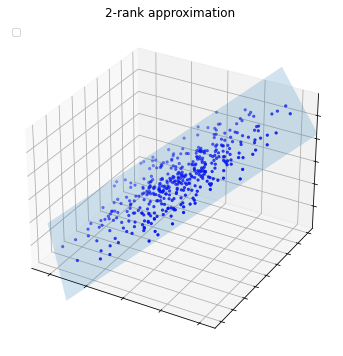}
\caption{A dataset in $\R^3$ (up) and its $1$-rank approximation (down-left) and its $2$-rank approximation (down-right) using SVD.}\label{fig:pca3Dapp}
\end{figure}

Up to here we have found all the principal components in order of significance.
The next step is to choose whether to keep all these components or discard those of lesser significance
and form with the remaining ones a new data in a lower dimensional subspace. Assume that we decide to keep only $k<d$ principal components.
Our new data matrix will be
$$
X_k = U_k\Sigma_k V_k^T = \sum_{j=1}^{k}\sigma_j u_{\cdot j}v_{\cdot j}^T
$$
which is a $k$-rank approximation of the original data matrix $X$, and is called the {\em feather matrix}. In fact, the new dataset is located on
a $k$-dimensional subspace of $\R^d$ spanned by orthonormal vectors $\{u_{\cdot1}, \ldots,u_{\cdot k}\}$.
This is a {\bf dimensionality reduction} algorithm based on SVD. See Figures \ref{fig:pca2Dapp} and \ref{fig:pca3Dapp} for illustrations in $2$ and $3$ dimensions, respectively .

\subsubsection*{Keywords and key sentences extraction}
Development of automatic procedures for {\em text summarization} is important
because people are overwhelmed by the tremendous amount of online information
and documents.

The goal of a text summarization algorithm is to {\em extract content
from a text document and present the most important content to the user in a
condensed form.}

We follow \cite{Elden:2019} and briefly discuss an algorithm based on SVD for automatically extracting key words and key sentences from a text.
Assume that we are given a text document, for example an article or a chapter of a book. Assume that the text contains
$m$ words and $n$ sentences\footnote{
Usually a couple of preprocessing steps should perform before summarization which are known as {\em stemming} and {\em stop words removing}.
In the stemming step the same word stem with different endings is represented by one token only.
For example words $\{$computation, computations, compute, computable, computing, computational$\}$ are represented by stem ``comput''.
Stop words such as $\{$a, an, and, or, if, then, any, as, about, able, \ldots$\}$ occur frequently in all texts and do not distinguish between different sentences. They should be
removed. We assume that special symbols, e.g., mathematics symbols are also removed.
}.
Normally $m$ is much larger than $n$. We form a $m\times n$ {\em term-sentence} matrix $A$ where its entry $a_{ij}$
is defined as the frequency of term (word) $i$
in sentence $j$.
Let us give the term $i$ the nonnegative {\em score} $u_i$ and the sentence $j$ the nonnegative {\em score} $v_j$.
The assignment of scores is made based on the following mutual reinforcement
principle:
\begin{quote}
A term should have a high score if it appears in many sentences
with high scores. A sentence should have a high score
if it contains many words with high scores.
\end{quote}
Using this criterion the score of term $i$ is proportional to the sum
of the scores of the sentences, weighted by frequencies, where it appears,
$$
\sigma u_i = \sum_{j=1}^{n}a_{ij} v_j, \quad i = 1,2,\ldots,m,
$$
where $\sigma>0$ is the proportional constant. Similarly, the score of sentence $j$ is proportional to the sum of scores
of its words weighted by their frequencies in the sentence,
$$
\sigma v_j = \sum_{i=1}^{m}a_{ij} u_i, \quad j = 1,2,\ldots,n.
$$
In matrix forms we have
$$
Av = \sigma u, \quad A^Tu = \sigma v
$$
for $u=[u_1,\ldots,u_m]^T$ (the score vector of words) and $v = [v_1,\ldots,v_n]^T$ (the score vector of sentences).
Collecting both equations we have
$$
AA^T u = \sigma^2 u, \quad A^TA v = \sigma^2v
$$
which show that $u$ is an eigenvector of $AA^T$ and $v$ is an eigenvector of $A^TA$ both corresponding to eigenvalue $\sigma^2$.
This means that $u$ is left singular vector and $v$ is a right singular value of $A$ corresponding to the same singular value $\sigma$.
If we choose the largest singular value, then we are guaranteed that the components
of $u$ and $v$ are nonnegative because the matrix $A$ has nonnegative entries.
Consequently, a simple algorithm based on SVD consists of the following steps:
\begin{itemize}
\item {\bf Step 1}: Apply a stemming and a stoping word algorithm on the text,
\item {\bf Step 2}: Construct the term-sentence matrix $A$,
\item {\bf Step 3}: Compute the SVD of $A$ and report the leading vectors $u$ and $v$ (first columns of $U$ and $V$, respectively) as word and sentence scores.
\end{itemize}
\noindent
The full SVD is not necessary for this algorithm as the first singular vectors are only required. Moreover, the term-sentence matrix is
sparse, so a SVD function for sparse matrices is numerically more efficient.
\vsp 

\subsection*{Classification of handwritten digits}

Classification of handwritten digits is a standard problem in pattern
recognition. A typical application is the automatic reading of zip codes on envelopes.
Here we follow the presentation of \cite{Elden:2019}
and give a classification technique using SVD.

The problem is as follows:
{\em Given a set of manually classified digits as a training set, classify a set
of unknown digits as the test set}.

In Figure \ref{fig:handwritten0} a sample of handwritten digits $0,1,\ldots,9$ is illustrated.
The images are downloaded from\footnote{\tt https://www.nist.gov/itl/products-and-services/emnist-dataset}.
Each image is a $28 \times 28$ grayscale image but we stack the columns of the image above each other so that each image is represented by a vector of size $784$. 

\begin{center}
  \includegraphics[scale=.8]{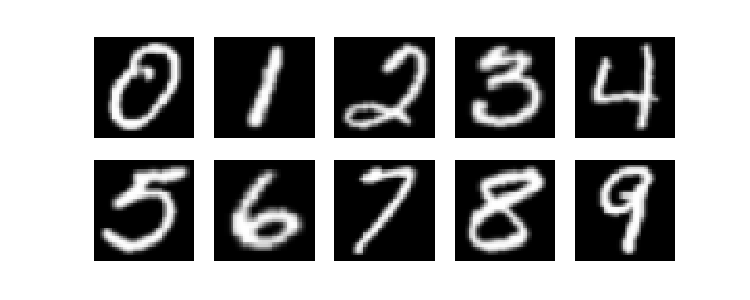}
  \captionof{figure}{A sample of handwritten digits $0,1,\ldots,9$.}\label{fig:handwritten0}
\end{center}

The training set of each kind of digits $0,1,\ldots,9$ can be considered as a matrix
$A$ of size $m\times n$ with $m = 784$ and $n =$ `the number of training digits of one kind'. Usually, $n\geqslant m$ although the case $n<m$ is also possible. So, we have a total of 10 training matrices, each corresponding to one of the digits $0, 1, \ldots, 9$.  
Assume that $A$ is one of them. The columns of $A$ span a linear subspace
of $\R^{784}$. However, this subspace cannot be expected to have a large dimension,
because columns of $A$ represent the same digit with different handwritings.
One can use SVD $A=U\Sigma V^T$ to
obtain an orthogonal basis for the column space of $A$ (or $\range(A)$) via the columns of factor $U$. Then any test image
of that kind can be well represented in terms of the orthogonal basis $\{u_{\cdot 1},u_{\cdot 2},\ldots,u_{\cdot m}\}$. Let us describe it in more detail.
We have
$$
A = U\Sigma V^T = \sum_{j=1}^{m}\sigma_j u_{\cdot j} v_{\cdot j}^T,
$$
thus the $\ell$-th column of $A$ (the $\ell$-th training digit) can be represented as
$$
a_{\cdot \ell} = \sum_{j=1}^{m} (\sigma_j v_{j\ell }) u_{\cdot j}
$$
which means that the coordinates of image $a_{\cdot \ell}$ in terms of orthogonal basis $\{u_{\cdot 1},u_{\cdot 2},\ldots,u_{\cdot m}\}$ are
$\sigma_j v_{j\ell }=:x_j$, i.e.
$$
a_{\cdot \ell}  = \sum_{j=1}^{m} x_j u_{\cdot j} = Ux.
$$
This means, in another point of view, that solving the least squares problem
$$
\min_{x}\|Ux - a_{\cdot \ell}\|_2
$$
results in an optimal vector $x=[\sigma_1 v_{1\ell },\ldots,\sigma_m v_{m\ell }]^T$ with residual $0$.
As we pointed out, most columns of $A$ are nearly linearly dependent as they represent different handwritten for the same digit. If we translate it to columns of $U$, this means that a few leading columns of $U$ should well capture the subspace. For this reason the orthogonal vectors $u_{\cdot 1},u_{\cdot 2},\ldots,$ are called
{\em singular images}. For example, in Figure \ref{fig:handwritten1}
the singular values and the first three singular images for
the training set $3$'s (containing the first $200$ images of the training set) are illustrated.
\begin{center}
  \includegraphics[scale=.55]{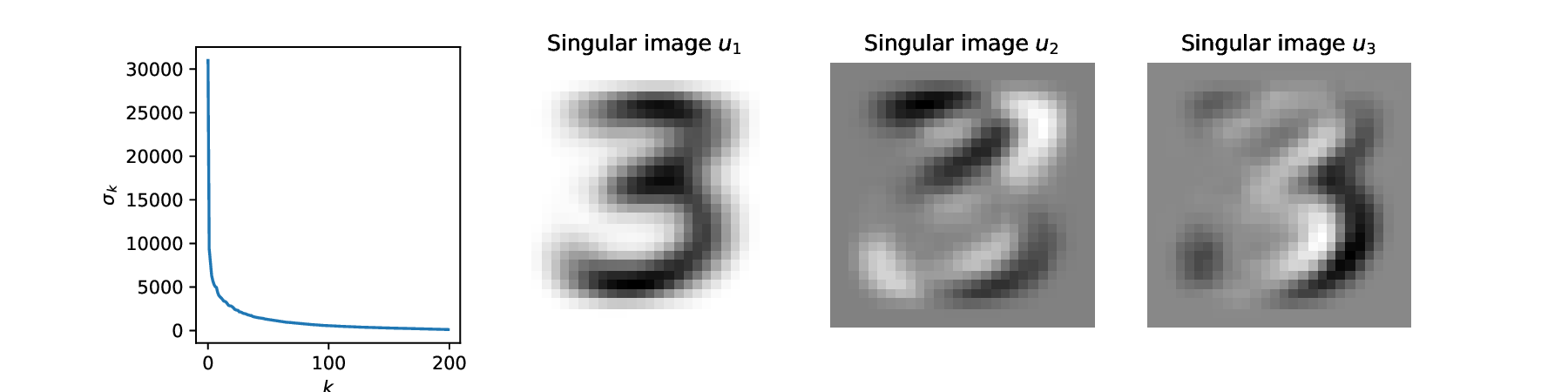}
  \captionof{figure}{Rapid decay of singular values, and the first three singular images.}\label{fig:handwritten1}
\end{center}
We observe that the first
singular image looks like a $3$, and the following singular images represent
the dominating variations of the training set around the first singular image.
Consequently, we can use
a $k$-rank approximation of $A$ where $k\ll m$. In this case each $a_{\cdot \ell}$ can be well represented by orthogonal basis
$\{u_{\cdot 1},u_{\cdot 2},\ldots,u_{\cdot k}\}$ in a least squares sense with a small residual.

Now, let $d\in\R^m$, $m=784$, be a digit  outside the training set (which is called a {\em test digit}). It is reasonable to assume that
$d$ can be well represented in terms of singular images $\{u_{\cdot 1},u_{\cdot 2},\ldots,u_{\cdot k}\}$ of its own type with a relatively small residual.
We must compute how well $d$ can be represented in
the $10$ different bases, each corresponding to one of the digits $0, 1, \ldots, 9$. This can be done by computing the residual vector in the least
squares problems of the type
\begin{equation*}\label{lstq0}
  \min_{x} \|U_kx-d\|_2
\end{equation*}
where $d$ represents a test digit and $U_k$ is a ${m\times k}$ matrix consisting of the first $k$ columns of $U$.
Since the columns of $U_k$ are orthonormal, the solution of this problem is $x = U_k^Td$ and the residual is
\begin{equation}\label{handw_res}
\|(I-U_k^TU_k)d\|_2.
\end{equation}
Altogether, a SVD-based classification algorithm of handwritten digits consists of two steps:
\begin{itemize}
\item {\bf Step 1 (training):} For the training set of known digits, compute the SVD of each set of digits of one kind (Compute ten SVDs for digits $0,1,\ldots,9$).
\item {\bf Step 2 (classification):} For a given test digit $d$, compute its relative residual in all ten bases using equation \eqref{handw_res}.
If a residual in a class is smaller than all the others, classify $d$ in that class.
\end{itemize}
\noindent
To test the algorithm, we download the {\em training} set consisting of $240000$ images ($24000$ images per any digit), and the {\em test} set consisting of $40000$ images ($4000$ images per digit) from the above-mentioned website. The document also contains two {\em label} files which are manual classifications of training and test sets. The training labels will be used to learn and the test labels to verify the algorithm. 

For each digit we consider using only the initial $200$ images out of the total $24000$ training images (ten matrices each of size $784\times 200$). However, we test the algorithm on all $40000$ test images and compare our classifications with the exact test labels. 
The percentages of the success of this SVD-based algorithm for each digit with different values $k = 3,6,\ldots,10$ (number of basis functions) are plotted in Figure \ref{fig:handwritten3}. It is left to readers to discuss their insights based on the observations from the figure. 

\begin{center}
  \includegraphics[scale=1]{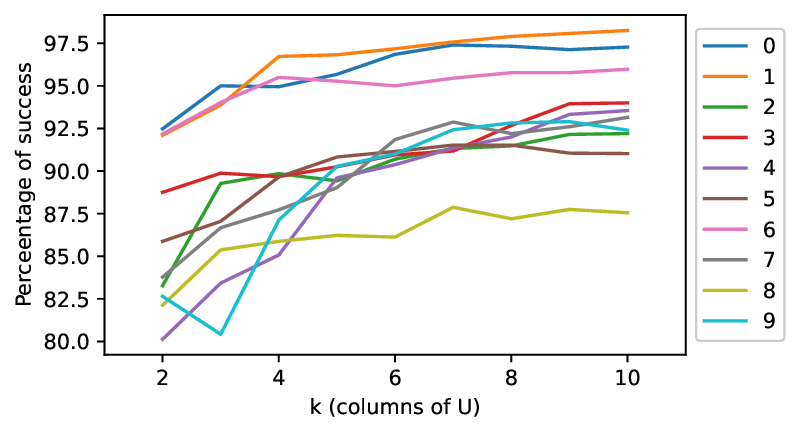}
  \captionof{figure}{Percentages of success ($y$-axis) in terms of number of basis functions ($x$-axis) for different digits $0,1,\ldots,9$.}\label{fig:handwritten3}
\end{center}

\begin{labexercise}
Implement your own Python code for the above classification algorithm.
Test your algorithm for different number of basis functions. 
Note that in the statement above, the equations you need are formulated for a single test vector $d$ (one ``flattened'' test image). However, if you use that approach in your code, it will take a very long time to run. In contrast, matrix-matrix multiplications are much more efficient on modern computers. To take full advantage of that, simply replace vector $d$ in equation \eqref{handw_res} with the entire test matrix. 
\end{labexercise}
\vsp

%% file: workouts_2.tex
\section*{Additional workouts}

\begin{workout}[Uniqueness of reduced QR factorization]
Assume that $A\in \R^{m\times n}$ has a
full rank. Prove that $A$ possesses a unique reduced QR factorization $A= Q_1R_1$ with the diagonal
entries of $R_1$ positive.
\end{workout}

\begin{workout}
Assume that $A\in \R^{m\times n}$ has rank $r <\min\{n, m\}$. Show that
for every $\epsilon > 0$ (no matter how much small), there exists a full-rank matrix $A_\epsilon\in \R^{m\times n}$
such that $\|A-A_\epsilon\|_2\leqslant \epsilon$.
\end{workout}

\begin{workout}
Let $A\in\R^{n\times n}$ be an arbitrary matrix. Find the closet orthogonal matrix $Q\in\R^{n\times n}$ to $A$ in Frobenius norm, i.e. solve the following problem:
$$
\min_{Q^TQ=I}\|A-Q\|_F.
$$
Hint: use SVD.
\end{workout}

\begin{workout}[\cite{Trefethen-Bau:1997}]
Let $P\in\R^{n\times n}$ be a nonzero projector. Show that $\|P\|_2\geqslant 1$ with equality if and only if $P$ is an orthogonal projector. 
\end{workout}

\begin{workout}[\cite{Trefethen-Bau:1997}]
Suppose that a square matrix $A$ has SVD $A=U\Sigma V^T$. Find the eigendecomposition of the symmetric matrix 
$$
B = \begin{bmatrix}
0 & A^T \\ A & 0
\end{bmatrix}.
$$
\end{workout}

\begin{workout}[\cite{Trefethen-Bau:1997}]
Suppose that $A$ is a $m\times n$ matrix with $m>n$ of the form 
$$
A = \begin{bmatrix}
B \\ C
\end{bmatrix}
$$
where $B$ is a nonsingular matrix of size $n\times n$ and $C$ is an arbitrary matrix of size $(m-n)\times n$. Prove that 
$
\|A^+\|_2 \leqslant \|B^{-1}\|_2.
$
\end{workout}

\begin{workout}[\cite{Datta:2010}]
Let $A$ be an $m\times n$ matrix of full rank $r = \min\{m,n\}$. If $B$ is another $m \times n$ matrix such
that $\|B- A\|_2 < \sigma_r$ then show $B$ has also full rank.
\end{workout}

\begin{workout}[\cite{Datta:2010}]
Show that the relative distance of a nonsingular matrix $A$
to the nearest singular matrix $B$ is
$$
\frac{\|A-B\|_2}{\|A\|_2}= \frac{1}{\cond_2(A)}.
$$
\end{workout}
\vsp 

%